\documentclass{amsart}

\usepackage{latexsym}
\usepackage{amsfonts,amsmath,amssymb,stmaryrd,amscd,amsxtra,amsthm,verbatim,calc}
\usepackage[all]{xy}

\RequirePackage[dvipsnames,usenames]{color}

\theoremstyle{plain}
\newtheorem{thm}{\bf Theorem}[section]
\newtheorem{pro}[thm]{\bf Proposition}
\newtheorem{lem}[thm]{\bf Lemma}
\newtheorem{defn}[thm]{\bf Definition}
\newtheorem{cor}[thm]{\bf Corollary}
\newtheorem{rem}[thm]{\bf Remark}

\newtheorem{ex}[thm]{\bf Example}

\newtheorem{defn-pro}[thm]{\bf Definition-Proposition}
\newtheorem{discussion}[thm]{\bf Discussion}

\DeclareMathOperator{\Hom}{Hom}
\DeclareMathOperator{\Ext}{Ext}
\DeclareMathOperator{\Ann}{Ann}

\DeclareMathOperator{\Proj}{Proj}

\DeclareMathOperator{\im}{Im}
\DeclareMathOperator{\Spec}{Spec}
\DeclareMathOperator{\cE}{\mathcal{E}}

\DeclareMathOperator{\cM}{\mathcal{M}}
\DeclareMathOperator{\cO}{\mathcal{O}}

\DeclareMathOperator{\cA}{\mathcal{A}}
\DeclareMathOperator{\cB}{\mathcal{B}}
\DeclareMathOperator{\cC}{\mathcal{C}}

\DeclareMathOperator{\Char}{char}
\DeclareMathOperator{\D}{D}
\DeclareMathOperator{\E}{E}

\DeclareMathOperator{\fm}{\mathfrak{m}}

\DeclareMathOperator{\dualx}{\omega^{\bullet}_X}

\DeclareMathOperator{\ch}{char}
\DeclareMathOperator{\dualy}{\omega^{\bullet}_Y}
\DeclareMathOperator{\dualz}{\omega^{\bullet}_Z}
\DeclareMathOperator{\dualp}{\omega^{\bullet}_{\mathbb{P}^n_k}}

\DeclareMathOperator{\Int}{Int}

\begin{document}
\title[Lyubeznik numbers of projective schemes]{Lyubeznik numbers of projective schemes}
\author[Wenliang Zhang]{Wenliang Zhang}
\address{Department of Mathematics, University of Michigan, Ann Arbor, MI 48109}
\email{wlzhang@umich.edu}
\subjclass[2000]{13D45, 14B15}
\keywords{local cohomology, Frobenius endomorphism, Lyubeznik number, projective scheme}

\begin{abstract}
Let $X$ be a projective scheme over a field $k$ and let $A$ be the local ring at the vertex of the affine cone of $X$ under some embedding $X\hookrightarrow\mathbb{P}^n_k$. We prove that, when $\ch(k)>0$, the Lyubeznik numbers $\lambda_{i,j}(A)$ are intrinsic numerical invariants of $X$, i.e., $\lambda_{i,j}(A)$ depend only on $X$, but not on the embedding.
\end{abstract}
\maketitle
\numberwithin{equation}{thm}
\section{Introduction}
Let $A$ be a local ring that contains a field and admits a surjection from an $n$-dimensional regular
local ring $(R,\mathfrak{m})$ containing a field. Let $I\subset R$ be the kernel of the surjection, and $k=R/{\mathfrak{m}}$ the residue field of $R$.  Lyubeznik numbers $\lambda_{i,j}(A)$ (Definition 4.1 in \cite{l1})
are defined to be $\dim_k(\Ext_R^i(k,H^{n-j}_I(R)))$. And it was proven in \cite{l1} that they are all finite and depend
only on $A,i$ and $j$, but neither on $R$, nor on the surjection $R\rightarrow A$. Lyubeznik numbers have been studied by a number of authors; see, for example, \cite{bb}, \cite{ls}, \cite{k1}, \cite{k2},  \cite{l3}, \cite{w}, and \cite{zh}.\par

Among all the things that make Lyubeznik numbers interesting is their close connection to topology. \par
Let $X$ be a complex variety with an isolated singular point $a$ ($X$ is smooth everywhere else). Let $A$ be the local ring of $X$ at $a$. Then it is proved in \cite{ls} that 
\[\lambda_{0,j}(A)=\dim_{\mathbb{C}}(H^j_{\{a\}}(X;\mathbb{C}))\]
and all other $\lambda_{i,j}(A)$ can be written in terms of $\lambda_{0,j}(A)$, where $H^j_{\{a\}}(X;\mathbb{C})$ denotes the $j$-th singular cohomology of $X$ with coefficients in $\mathbb{C}$ supported  at $\{a\}$. This result is quite striking. It indicates that there is some connection between $\lambda_{i,j}(A)$ and the {\itshape singular} topology of a neighborhood of the point $a$; however, $\lambda_{i,j}(A)$ are defined purely algebraically. \par 
Given a $d$-dimensional projective scheme $X$ over a field $k$ of any characteristic (for simplicity assuming $k$ is algebraically closed), under an embedding $X\hookrightarrow \mathbb{P}^n_k$, one can write $X=\Proj(k[x_0,\dots,x_n]/I)$, where $I$ is a homogeneous ideal of $R=k[x_0,\dots,x_n]$. Let $A=(R/I)_{(x_0,\dots,x_n)}$, i.e., the local ring at the vertex of the affine cone of $X$ under $X\hookrightarrow \mathbb{P}^n_k$. Let $X_1,\dots,X_s$ be the $d$-dimensional irreducible components of $X$ and let $\Gamma_X$ be the graph on the vertices $X_1,\dots,X_s$ in which $X_i$ and $X_j$ are joined by an edge if and only if $\dim(X_i\cap X_j)=d-1$. Then \cite[Theorem 2.7]{zh} shows that $\lambda_{d+1,d+1}(A)$ equals the number of connected components of $\Gamma_X$. This result is somewhat intriguing in that, not only does it indicate that some topological information of $X$ is encoded in $\lambda_{d+1,d+1}(A)$, but equally importantly it shows that $\lambda_{d+1,d+1}(A)$, whose definition {\itshape a priori} requires an embedding $X\hookrightarrow\mathbb{P}^n_k$, depends only on $X$ but not on the embedding. \par

The goal of this paper is to establish the embedding-independence of all $\lambda_{i,j}(A)$ in positive characteristic.

\begin{thm}[Main Theorem]
\label{main-thm}
Let $X$ be a projective scheme over a field $k$ of positive characteristic, and let $A$ be the local ring at the vertex of the affine cone over $X$ for some embedding of $X$ into a projective space. Then $\lambda_{i,j}(A)$ depends only on $X$, $i$ and $j$, but not on the embedding.
\end{thm}

This gives a positive answer to the open problem posed in \cite[p. 133]{l2} in the case when $\ch(k)>0$. To the best of our knowledge, in characteristic 0 the problem is still open.\par

Let $R$ be a noetherian ring of positive characteristic $p$ and let $k$ be a coefficient field of $R$ which for the purpose of this introduction we assume to be perfect. Let $M$ be an $R$-module equipped with a $p$-linear structure $f$, i.e., a map of abelian groups $f:M\to M$ such that $f(rm)=r^pf(m)$ for all $r\in R,m\in M$. The stable part of $M$, denoted by $M_s$, is defined to be
$$\bigcap ^{\infty}_{i=1}f^i(M)$$
(in general, this is not an $R$-module, but a $k$-vector-space). The operation of `taking the stable part' has played an important role in the study of local cohomology in positive characteristic. For example, it has been used to estimate local cohomological dimensions (cf. \cite{lcd}) and to prove vanishing results of local cohomology (cf. \cite{l5}), etc. In this paper, we will give a description of $\lambda_{i,j}(A)$ ($A$ as in Theorem \ref{main-thm}) in terms of the stable part of a certain module $\cE^{i,j}(R/I)$ (defined in the next paragraph).\par

Let $\eta:X\hookrightarrow \mathbb{P}^n_k$ be an embedding of $X$.  Let $R=k[x_0,\cdots,x_n]$ and $I$ be the defining ideal of $X$ ($I$ is a homogeneous ideal of $R$). Set $$\cE^{i,j}(M)=\Ext^{n+1-i}_R(\Ext^{n+1-j}_R(M,\Omega),\Omega)$$
where $\Omega=R(-n-1)$ is the graded canonical $R$-module, for all graded $R$-modules $M$. The module $\cE^{i,j}(M)$ is naturally graded and the degree-$l$ piece of $\cE^{i,j}(M)$ is a finite-dimensional $k$-vector space for every integer $l$.\par

The $p$-linear structure $\varphi$ on $\cE^{i,j}(R/I)$ induced by the map $R/I\xrightarrow{\hat{r}\mapsto \hat{r}^p}R/I$ satisfies $\deg(\varphi(m))=p\deg(m)$ (see Section 5 for details) and hence induces a $p$-linear structure $\varphi_0$ on $\cE^{i,j}(R/I)_0$, the degree-0 piece of $\cE^{i,j}(R/I)$. The connection between $\lambda_{i,j}(A)$  ($A$ as in Theorem \ref{main-thm}) and $\cE^{i,j}(R/I)_0$ lies in the following theorem.
\begin{thm}
\label{stable-M-i-j}
Let $R,I,\cE^{i,j}(R/I)$ be as above and $A=(R/I)_{(x_0,\cdots,x_n)}$ be the local ring at the vertex of the affine cone over $X$ under an embedding $\eta:X\hookrightarrow \mathbb{P}^n_k$, where $k$ is a perfect field of positive characteristic $p$. Then
$$\lambda_{i,j}(A)=\dim_k((\cE^{i,j}(R/I)_0)_s).$$
\end{thm}

Given Theorem \ref{stable-M-i-j}, it should be clear that, to establish Theorem \ref{main-thm}, it suffices to prove that following two theorems. (Note that in Theorem \ref{indp-M-i-j}, we do not assume that $k$ is of positive characteristic.)

\begin{thm}[Independence of $\cE^{i,j}(R/I)_0$]
\label{indp-M-i-j}
Let $X,R,I,\cE^{i,j}(R/I)$ be as above and assume that $k$ is a field of any characteristic. Then $\cE^{i,j}(R/I)_0$ depends only on $X$, $i$ and $j$, but not on the embedding $\eta:X\hookrightarrow \mathbb{P}^n_k$.
\end{thm}

\begin{thm}[Independence of $\varphi_0$]
\label{indp-f-M}
The $p$-linear structure $\varphi_0$ on $\cE^{i,j}(R/I)_0$ is independent of the embedding $\eta:X\hookrightarrow \mathbb{P}^n_k$.
\end{thm}

The paper is organized as follows. In Section 2, we collect some basic facts that will be used in subsequent sections. In Section 3,
we introduce a covariant functor $\E^{s,t}$ from the category of quasi-coherent sheaves on a scheme $Z$ to the category of
$\Gamma(Z,\mathcal{O}_Z)$-modules and prove that this functor commutes with finite morphisms (Definition \ref{E-ij} and Corollary
\ref{f-finite-gamma}). In Section 4, we study $\cE^{i,j}(-)$ in detail, in particular, we will prove Theorem \ref{indp-M-i-j} as a corollary to Proposition \ref{M-ij}. All results in Section 2, 3 and 4 are characteristic-free. Results in the rest of this paper are in positive characteristic. In Section 5, we develop a graded theory of the Frobenius endomorphism and functors associated with it. Section 6 and 7 are devoted to proving Theorem \ref{indp-f-M}. Finally in Section 8, we establish Theorem \ref{stable-M-i-j}. \par

\textbf{Notations}. Throughout this paper, all schemes will be assumed to be noetherian and we will follow notations as in \cite{rd}; in particular, the usual sheaf Hom (or sheaf Ext) will be denoted by $\underline{\Hom}$ (or by $\underline{\Ext}$).\par
Let $Y$ be any scheme and let $\mathcal{F}$ be a quasi-coherent sheaf of $\mathcal{O}_Y$-modules. Let $\mathcal{C}^{\bullet}$ be a complex of quasi-coherent sheaves of $\mathcal{O}_Y$-modules and let $\mathcal{I}^{\bullet}$ be an injective resolution of $\mathcal{C}^{\bullet}$, i.e., a complex of injectives with a quasi-isomorphism $\mathcal{C}^{\bullet}\to \mathcal{I}^{\bullet}$. Then by $\underline{\Ext}^i_Y(\mathcal{F},\mathcal{C}^{\bullet})$ ($\Ext^i_Y(\mathcal{F},\mathcal{C}^{\bullet})$ respectively), we mean the $i$-th cohomology of the complex $\underline{\Hom}_{Y}(\mathcal{F},\mathcal{I}^{\bullet})$ ($\Hom_{Y}(\mathcal{F},\mathcal{I}^{\bullet})$ respectively). We omit the subscript $Y$ if $Y$ is clear from the context.\par
For any integer $t$, we denote by $\mathcal{C}^{\bullet}[t]$ the complex such that $(\mathcal{C}^{\bullet}[t])^l=\mathcal{C}^{l+t}$ for all integers $l$ (i.e., $\mathcal{C}^{\bullet}[t]$ is the complex $\mathcal{C}^{\bullet}$ shifted to the left by $t$ places). And $\mathcal{F}[t]$, for any sheaf $\mathcal{F}$, denotes the complex with $\mathcal{F}$ concentrated in homological degree $-t$, i.e., $(\mathcal{F}[t])^l=\begin{cases} \mathcal{F} & l=-t\\ 0 & \text{otherwise}\end{cases}$.\par

For any scheme $Z$, we will use $\Gamma(Z)$ to denote $\Gamma(Z,\cO_Z)$, the ring of global sections on $Z$. For any scheme morphism $(f,f^{\#}):X\to Y$, we will use 
$$f^{\Gamma}_*:\Gamma(X)\text{-mod}\to \Gamma(Y)\text{-mod}$$
to denote the restriction functor from the category of $\Gamma(X)$-modules to the category of $\Gamma(Y)$-modules, i.e., a $\Gamma(X)$-module is regarded as a $\Gamma(Y)$-module via $f^{\#}:\Gamma(Y)\to \Gamma(X)$.\par

Let $C$ be an arbitrary commutative noetherian ring. Then there is an equivalence between the category of quasi-coherent sheaves on $\Spec(C)$ and the category of $C$-modules; for example, the sheaf $\underline{\Ext}^l_{\Spec(C)}(\widetilde{M},\widetilde{N})$ corresponds under this equivalence to $\Ext^l_C(M,N)$ for any $C$-modules $M$ and $N$. Hence we will treat $\widetilde{M}$ and $M$ as the same object for all $C$-modules $M$, assuming no confusion will arise.\par

For any commutative noetherian ring $C$ of positive characteristic $p$, an integer $e\geq 1$ and an ideal $J$ in $C$, we denote by $J^{[p^e]}$ the ideal of $R$ generated by $\{r^{p^e}|r\in J\}$.\par

When $C$ is a graded commutative ring and $M$ is a graded $C$-module, for any integer $t$, we denote by $M(t)$ the graded $C$-module whose degree-$l$ piece is $M_{l+t}$, i.e., $M$ degree-shifted downward by $t$.\par 
When $R=k[x_0,\dots,x_n]$ is considered as a graded ring, the grading is always the standard one, i.e., $\deg(x_i)=1$ for $0\leq i\leq n$ and $\deg(c)=0$ for all $c\in k$. The irrelevant maximal ideal $(x_0,\dots,x_n)$ will be denoted by $\fm$. And we will use $\Omega$ to denote $R(-n-1)$, the graded canonical module of $R$. For each graded $R$-module $M$, we set
\[\cE^{i,j}(M)=\Ext^{n+1-i}_R(\Ext^{n+1-j}_R(M,\Omega),\Omega).\]

\section{Preliminaries}
In this section, we collect for future references some basic facts which will be used in subsequent sections.\par

Let $S$ be a graded noetherian ring and let $M$ be a graded $S$-module. In the following remark $M_{(f)}$ denotes the homogeneous localization of a graded module $M$ with respect to the multiplicative system $\{1,f,f^2,\dots\}$ (i.e., the degree-0 part of $M_f$), where $f$ is a homogeneous element; $\widetilde{M}$ denotes the sheaf on $\Proj(S)$ associated to $M$; and,
$$\sideset{^*}{_S}\Hom(M,N):=\oplus_n\Hom_S(M,N)_n,$$
where $\Hom_S(M,N)_n$ is the set of homomorphisms of degree $n$ (see \S2 in \cite{ega2} for details; in \cite{ega2} $\sideset{^*}{}\Hom _S(M,N)$ is denoted simply by $\Hom_S(M,N)$). If $M$ is finitely generated, $\sideset{^*}{}\Hom _S(M,N)$ coincides with $\Hom_S(M,N)$ in the usual sense.

\begin{rem}[(2.5.12 in \cite{ega2})]
\label{ega2}
{\rm Let $S$ be a graded noetherian ring, $M$ and $N$ two graded $S$-modules, and $f\in S_d$ ($d>0$).\par
One can define a canonical functorial $S_{(f)}$-module homomorphism
$$\mu_f: (\sideset{^*}{_S}\Hom(M,N))_{(f)}\rightarrow \Hom_{S_{(f)}}(M_{(f)},N_{(f)}) $$
by sending $u/f^l$, where $u$ is a homomorphism of degree $ld$ to the homomorphism $M_{(f)}\rightarrow N_{(f)}$ which maps
$x/f^m$ ($x\in M_{md}$) to $u(x)/f^{l+m}$.\par
For $g\in S_e$ ($e>0$), moreover, one has a commutative diagram
$$\begin{CD}
(\sideset{^*}{}\Hom_S(M,N))_{(f)} @>\mu_f>>        \Hom_{S_{(f)}}(M_{(f)},N_{(f)})  \\
      @VVV                                     @VVV                  \\
(\sideset{^*}{}\Hom_S(M,N))_{(fg)} @>\mu_{fg}>>    \Hom_{S_{(fg)}}(M_{(fg)},N_{(fg)})
\end{CD}$$
The vertical arrows are the canonical homomorphisms.\par
Furthermore, these $\mu_f$ define a canonical functorial homomorphism of $\mathcal{O}_X$-modules
$$\mu: \widetilde{\sideset{^*}{_S}\Hom(M,N)}\rightarrow \underline{\Hom}_{\mathcal{O}_X}(\widetilde{M},\widetilde{N}),$$
where $X=\Proj(S)$.}
\end{rem}

We denote by $S_+$ the ideal of $S$ generated by the elements of positive degrees.
\begin{pro}[Proposition 2.5.13 in \cite{ega2}]
 Suppose, for a graded noetherian ring $S$, the ideal $S_+$ is generated by $S_1$. Then
$$\mu: \widetilde{\sideset{^*}{_S}\Hom(M,N)}\rightarrow \underline{\Hom}_{\mathcal{O}_X}(\widetilde{M},\widetilde{N})$$
is an isomorphism when $M$ is finitely generated and $X=\Proj(S)$.
\end{pro}

The following proposition should be well-known to experts. Since we could not find a proof in the literature, we include a proof here.

\begin{pro}
\label{sheaf-ext}
Let $S$ be a graded Noetherian ring. Suppose that $S_+$ is generated by $S_1$. Let $M$ and $N$ be two graded $S$-modules. Then the homomorphism $\mu$ in Remark \ref{ega2} induces an isomorphism
$$\mu: \widetilde{\sideset{^*}{^i_S}\Ext(M,N)}\xrightarrow{\cong}\underline{\Ext}^i_{\mathcal{O}_X}(\widetilde{M},\widetilde{N}),$$
when $M$ is finitely generated and $X=\Proj(S)$.
\end{pro}

\begin{proof}[Proof]
$M$ has a resolution with degree-preserving differentials $F^{\bullet}\to M\to 0$ in which all $F_i$ are finitely generated graded free $S$-modules. Since $\widetilde{\cdot}$ is an exact functor, we have an induced locally-free resolution $\widetilde{F}^{\bullet}\to \widetilde{M}\to 0$ for $\widetilde{M}$.
According to Proposition 6.5 in Chapter 3 in \cite{ag}, $\underline{\Ext}^i_{\mathcal{O}_X}(\widetilde{M},\widetilde{N})$ is the $i$-th cohomology sheaf of the induced complex
$\underline{\Hom}_{\cO_X}(\widetilde{F}^{\bullet},\widetilde{N})$. Applying $\sideset{^*}{_S}\Hom(\cdot, N)$ to the resolution of $M$, we have a complex ${^*\Hom}_S(F^{\bullet},N)$.
This induces a complex of sheaves $\widetilde{{^*\Hom}_S(F^{\bullet},N)}$. Since $\widetilde{\cdot}$ is exact,
$\widetilde{\sideset{^*}{^i_S}\Ext(M,N)}$ is the $i$-th cohomology sheaf of the complex $\widetilde{{^*\Hom}_S(F^{\bullet},N)}$. Since the isomorphism $\mu:\widetilde{\sideset{^*}{_S}\Hom(\cdot,N)}\to \underline{\Hom}_{\mathcal{O}_X}(\widetilde{\cdot},\widetilde{N})$ is functorial, we have an isomorphism of complexes
$$\widetilde{{^*\Hom}_S(F^{\bullet},N)}\xrightarrow{\sim}\underline{\Hom}_{\cO_X}(\widetilde{F}^{\bullet},\widetilde{N})$$ Therefore, the induced maps on homology are isomorphisms.
\end{proof}

Consider $R=k[x_0,\dots,x_n]$ as a graded ring with the standard grading. Let $N_1$ and $N_2$ be graded $R$-modules. The grading on $\sideset{^*}{}\Hom_R(N_1,N_2)$ is given as follows: a homogeneous element of degree $l$ is a homomorphism $\phi\in \Hom_R(N_1,N_2)$ such that $\deg(\phi(a))=\deg(a)+l$ for any homogeneous element $a\in N_1$. This induces a grading on $\sideset{^*}{^t_R}\Ext(N_1,N_2)$ for all integers $t$. When $N_1$ is a finitely generated graded $R$-module, one has $\sideset{^*}{}\Hom_R(N_1,N_2)=\Hom_R(N_1,N_2),$
and hence
$\sideset{^*}{^t_R}\Ext(N_1,N_2)=\Ext^t_R(N_1,N_2)$.
Therefore, when $N_1$ is finitely generated, we will not distinguish $\sideset{^*}{}\Ext^t_R(N_1,N_2)$ and $\Ext^t_R(N_1,N_2)$, and just write $\Ext^t_R(N_1,N_2)$ (with the same grading as on $\sideset{^*}{^t_R}\Ext(N_1,N_2)$ kept in mind). In particular, in what follows, we will write $$\Ext^{n+1-j}_R(R/I,\Omega)$$ and $$\Ext^{n+1-i}_R(\Ext^{n+1-j}_R(R/I,\Omega),\Omega)$$ with $^*$ dropped.\par
It is easy to see that a degree-preserving homomorphism $N_1\to N_2$ induces degree-preserving homomorphisms $\Ext^i_R(M,N_1)\to \Ext^i_R(M,N_2)$ (or \newline $\Ext^i_R(N_2,M)\to \Ext^i_R(N_1,M)$) for any graded $R$-module $M$.\par

Next, we collect some preliminary facts on dualizing complexes and Grothendieck Duality from \cite{rd}.\par
For any scheme $Y$, we will use $\mathcal{D}_c(Y)$ ($\mathcal{D}^+_c(Y)$, $\mathcal{D}^-_c(Y)$, respectively) to denote the the derived category of complexes (bounded-below complexes, bounded-up complexes, respectively) of sheaves $\cO_Y$-modules with coherent cohomology; similarly, we will use $\mathcal{D}_{qc}(Y)$ ($\mathcal{D}^+_{qc}(Y)$, $\mathcal{D}^-_{qc}(Y)$, respectively) to denote the the derived category of complexes (bounded-below complexes, bounded-up complexes, respectively) of sheaves $\cO_Y$-modules with quasi-coherent cohomology.\par
We recall the definition and some basic facts of dualizing complexes.\par

\begin{defn}[Definition on page 258 in \cite{rd}]
For any (locally) noetherian scheme $Y$, a complex $\mathcal{R}^{\bullet}\in \mathcal{D}^+_c(Y)$ with finite injective dimension is called a dualizing complex if the natural map
$$\cO_Y\to \underline{\underline{R}}\ \underline{\Hom}^{\bullet}_Y(\underline{\underline{R}}\ \underline{\Hom}^{\bullet}_Y(\cO_Y,\mathcal{R}^{\bullet}),\mathcal{R}^{\bullet})$$
is an isomorphism.
\end{defn}

The existence of dualizing complexes on a scheme $Y$ puts restrictions on $Y$, for instance, $Y$ must be universally catenary and have finite Krull dimension (\cite[page 300]{rd}). Sufficient conditions for the existence of dualizing complexes can also be found in \cite[page 299]{rd}, for example, any noetherian scheme of finite type over a field admits dualizing complexes.\par
 
If $f:X\to Y$ is a scheme morphism of finite type where both $X$ and $Y$ admit dualizing complexes, then (by \cite[Corollary 3.4 in Chapter VII]{rd}) there is a functor (called {\itshape twisted inverse image functor} by some authors)\[f^!:\mathcal{D}^+_{qc}(Y)\rightarrow \mathcal{D}^+_{qc}(X).\] 
If $\mathcal{R}^{\bullet}$ is a dualizing complex on $Y$, then $f^!\mathcal{R}^{\bullet}$ is also a dualizing complex on $X$ by \cite[page 299]{rd}.\par

If $f:X\to Y$ is a finite morphism of noetherian schemes, then $f^!$ can be defined explicitly as (\cite[Definition on page 165]{rd})
$$f^!(\cdot)=\bar{f}^*\underline{\underline{R}}\ \underline{\Hom}_{Y}(f_*\mathcal{O}_{X}, \cdot),$$
where $\bar{f}$ denotes the induced morphism
$$(X,\mathcal{O}_{X})\rightarrow (Y, f_*\mathcal{O}_{X}).$$

\begin{thm}[Duality for finite morphisms, Theorem 6.7 in Chapter 3 \cite{rd}]
\label{duality}
Let $f:X\to Y$ be a finite morphism of noetherian schemes of finite Krull dimension. Then the duality morphism
$$\theta_f:\underline{\underline{R}}f_*\underline{\underline{R}}\ \underline{\Hom}^{\bullet}_X(\mathcal{F}^{\bullet}, f^!\mathcal{G}^{\bullet})\to \underline{\underline{R}}\ \underline{\Hom}^{\bullet}_Y(\underline{\underline{R}}f_*\mathcal{F}^{\bullet}, \mathcal{G}^{\bullet})$$
is an isomorphism for $\mathcal{F}^{\bullet}\in \mathcal{D}^-_{qc}(X)$ and $\mathcal{G}^{\bullet}\in \mathcal{D}^+_{qc}(Y)$.
\end{thm}

\begin{rem}
\label{refined-finite-duality}
{\rm Let $\mathcal{F}^{\bullet}$ and $\mathcal{G}^{\bullet}$ be as in Theorem \ref{duality}. It is clear that $\underline{\underline{R}}f_*\mathcal{F}^{\bullet}=f_*\mathcal{F}^{\bullet}$, since $f$ is finite. Moreover, if $\mathcal{F}^{\bullet}\in\mathcal{D}^-_{c}(X)$ and $\mathcal{G}^{\bullet}\in\mathcal{D}^+_{c}(X)$, then according to \cite[Proposition 3.3 in Chapter II]{rd}, 
\[\underline{\underline{R}}\ \underline{\Hom}^{\bullet}_X(\mathcal{F}^{\bullet}, f^!\mathcal{G}^{\bullet})\in \mathcal{D}_{c}(X).\] 
And hence, when $\mathcal{F}^{\bullet}\in \mathcal{D}^-_{c}(X),\mathcal{G}^{\bullet}\in \mathcal{D}^+_{c}(X)$, we may replace $\underline{\underline{R}}f_*$ by $f_*$ in Theorem \ref{duality} since $f$ is finite, i.e., we have an isomorphism
\[f_*\underline{\underline{R}}\ \underline{\Hom}^{\bullet}_X(\mathcal{F}^{\bullet}, f^!\mathcal{G}^{\bullet})\to \underline{\underline{R}}\ \underline{\Hom}^{\bullet}_Y(f_*\mathcal{F}^{\bullet}, \mathcal{G}^{\bullet})\]
for $\mathcal{F}^{\bullet}\in \mathcal{D}^-_{c}(X)$ and $\mathcal{G}^{\bullet}\in \mathcal{D}^+_{c}(Y)$.}
\end{rem}

Let $X$ be a scheme, let $\mathcal{F}$ be a coherent sheaf of $\mathcal{O}_X$-modules and let $\mathcal{G}^{\bullet}\in \mathcal{D}^+_{qc}(X)$ be a complex of sheaves of $\mathcal{O}_X$-modules with quasi-coherent cohomology. Then $\Ext^i_X(\mathcal{F},\mathcal{G}^{\bullet})$ (defined as the $i$th cohomology of $\Hom^{\bullet}_X(\mathcal{F},\mathcal{I}^{\bullet})$ where $\mathcal{I}^{\bullet}$ is a complex of injective sheaves which is quasi-isomorphic to $\mathcal{G}^{\bullet}$) is naturally a $\Gamma(X)$-module.\par

\begin{cor}
\label{duality-sheaf-ext}
Let $f:X\to Y$ be a finite morphism of noetherian schemes of finite Krull dimension. Let $\mathcal{F}$ be a coherent sheaf of $\mathcal{O}_X$-modules and $\mathcal{G}^{\bullet}\in\mathcal{D}^+_{c}(Y)$. Then there is a natural isomorphism of sheaves of $\mathcal{O}_Y$-modules
\begin{equation}
\label{sheaf-ext-duality}
f_*\underline{\Ext}^i_X(\mathcal{F},f^!\mathcal{G}^{\bullet})\xrightarrow{\sim} \underline{\Ext}^i_Y(f_*\mathcal{F},\mathcal{G}^{\bullet}),
\end{equation}
and a natural isomorphism of $\Gamma(Y)$-modules
\begin{equation}
\label{ext-duality}
f^{\Gamma}_*\Ext^i_X(\mathcal{F},f^!\mathcal{G}^{\bullet})\xrightarrow{\sim}\Ext^i_Y(f_*\mathcal{F},\mathcal{G}^{\bullet}) .
\end{equation}
\end{cor}
\begin{proof}[Proof]
Let $\sigma:\mathcal{G}^{\bullet}\to \mathcal{I}^{\bullet}$ be an injective resolution (i.e., $\sigma$ is a quasi-isomorphism and $\mathcal{I}^{\bullet}$ is a complex of injectives), then
$$f^!\mathcal{G}^{\bullet}=\bar{f}^*\underline{\underline{R}}\ \underline{\Hom}^{\bullet}_Y(f_*\mathcal{O}_X,\mathcal{G}^{\bullet})=\bar{f}^*\underline{\underline{R}}\ \underline{\Hom}^{\bullet}_Y(f_*\mathcal{O}_X,\mathcal{I}^{\bullet})=f^!\mathcal{I}^{\bullet}.$$
Since $f^!$ is right-adjoint to an exact functor $f_*$ ($f$ is a finite morphism), according to \cite[Proposition 2.3.10]{weibel} $f^!$ sends injectives to injectives; hence $f^!\mathcal{I}^{\bullet}$ is a complex of injectives on $X$. Consider the following diagram of sheaves of $\cO_Y$-modules
$$\begin{CD}
\cdots @>>> f_*\underline{\Hom}_X(\mathcal{F},f^!\mathcal{I}^i) @>>> f_*\underline{\Hom}_X(\mathcal{F},f^!\mathcal{I}^{i+1}) @>>> \cdots \\
@.          @VV\sim V                    @VV\sim V       @.\\
\cdots @>>> \underline{\Hom}_Y(f_*\mathcal{F},\mathcal{I}^i) @>>> \underline{\Hom}_Y(f_*\mathcal{F},\mathcal{I}^{i+1}) @>>> \cdots
\end{CD}$$
where the vertical isomorphisms are given by Theorem \ref{duality} and Remark \ref{refined-finite-duality}, and the commutativity of the diagram follows from the functoriality of the isomorphism in Theorem \ref{duality}. Since $f_*\underline{\Ext}^i_X(\mathcal{F},f^!\mathcal{G}^{\bullet})$ is the $i$-th cohomology of the second row and $\underline{\Ext}^i_Y(f_*\mathcal{F},\mathcal{G}^{\bullet})$ is the $i$-th cohomology of the second row, this finishes the proof that  the map (\ref{sheaf-ext-duality}) exists and is an isomorphism.\par
The above diagram induces a commutative diagram of $\Gamma(Y)$-modules
$$\begin{CD}
\cdots @>>> \Gamma(Y,f_*\underline{\Hom}_X(\mathcal{F},f^!\mathcal{I}^i)) @>>> \Gamma(Y,f_*\underline{\Hom}_X(\mathcal{F},f^!\mathcal{I}^{i+1})) @>>> \cdots \\
@.          @VV\sim V                    @VV\sim V       @.\\
\cdots @>>> \Gamma(Y,\underline{\Hom}_Y(f_*\mathcal{F},\mathcal{I}^i)) @>>> \Gamma(Y,\underline{\Hom}_Y(f_*\mathcal{F},\mathcal{I}^{i+1})) @>>> \cdots
\end{CD}$$
Since the pre-sheaf $\underline{\Hom}$ is indeed a sheaf (page 109 in \cite{ag}) and $\Gamma(Y,f_*\mathcal{G})=f^{\Gamma}_*(\Gamma(X,\mathcal{G}))$ for every sheaf $\mathcal{G}$ on $X$, one has, in turn, the following commutative diagram
$$\begin{CD}
\cdots @>>> f^{\Gamma}_*\Hom_X(\mathcal{F},f^!\mathcal{I}^i) @>>> f^{\Gamma}_*\Hom_X(\mathcal{F},f^!\mathcal{I}^{i+1}) @>>> \cdots \\
@.          @VV\sim V                    @VV\sim V       @.\\
\cdots @>>> \Hom_Y(f_*\mathcal{F},\mathcal{I}^i) @>>> \Hom_Y(f_*\mathcal{F},\mathcal{I}^{i+1}) @>>> \cdots
\end{CD}$$
which completes the proof that the map (\ref{ext-duality}) exists and is an isomorphism.
\end{proof}

Next, we consider a graded version of inverse limits.\par
Let $S$ be a graded commutative ring, then one can define
$\sideset{^*}{}\varprojlim$ as a graded version of $\varprojlim$ of
an inverse system of graded $S$-modules as follows. Let
$\{M_i,\theta_{ji}:M_j\to M_i\}$ be an inverse system of graded
$S$-modules where $\theta_{ji}$ are degree-preserving $S$-module
homomorphisms. Define
$$(\sideset{^*}{}\varprojlim_iM_i)_n=\varprojlim_i(M_i)_n,$$ where
$\varprojlim_i(M_i)_n$ is the inverse limit of the inverse system
$\{(M_i)_n,(\theta_{ji})_n:(M_j)_n\to (M_i)_n\}$. Then define
$$\sideset{^*}{}\varprojlim_iM_i=\oplus_{n\in \mathbb{Z}}(\sideset{^*}{}\varprojlim_iM_i)_n.$$
Then $\sideset{^*}{}\varprojlim_iM_i$ is naturally a graded $S$-module.\par
In general, $\sideset{^*}{}\varprojlim$ and $\varprojlim$ are different as seen in the following example.

\begin{ex}
\label{example} Let $R=k[x_0,\dots,x_n]$ be a polynomial ring in
$n+1$ variables over a field $k$ with $\Char(k)=p>0$ and $J$ be a
proper homogeneous ideal of $R$. Then
$$\sideset{^*}{}\varprojlim_iR/J^{[p^i]}=R,\ {\rm but}\ \varprojlim_iR/J^{[p^i]}=\hat{R}^J,$$
where $\hat{R}^J$ is the $J$-adic completion of $R$.
\end{ex}
\begin{proof}[Proof]
The second assertion, $\varprojlim_iR/J^{[p^i]}=\hat{R}^J$, is clear.\par
To prove that $\sideset{^*}{}\varprojlim_iR/J^{[p^i]}=R$, it suffices to show
$$\varprojlim_i(R/J^{[p^i]})_n=R_n$$
for every integer $n\geq 0$.\par For each $n\geq 0$, there exists an
integer $i(n)$ such that $p^j>n$ when $j\geq i(n)$. Then
$(R/J^{[p^j]})_n=R_n$ for all $j\geq i(n)$, and hence the inverse system
$\{(R/J^{[p^i]})_n\}$ is indeed the following
\begin{equation}
\label{sn}
(R/J^{[p]})_n\leftarrow (R/J^{[p^2]})_n\leftarrow\cdots\leftarrow R_n\xleftarrow{=}R_n\xleftarrow{=}\cdots
\end{equation}
It is clear that the inverse limit of (\ref{sn}) is $R_n$.
\end{proof}
One may also define $\sideset{^*}{}\varinjlim$ as a graded version
of $\varinjlim$ similarly. However, by doing so, one does not get
anything new, since (it is easy to check that) $\varinjlim_iM_i$ is
naturally graded with respect to the grading on $M_i$'s whenever all
$M_i$ are graded and $\theta_{ij}$ are degree-preserving.\par

In general, this graded inverse limit commutes with ${^*\Hom}$ in the following sense.

\begin{pro}
\label{graded-inverse}
Let $S$ be a graded commutative ring. Let $N$ be a graded $S$-module and $\{M_i,\theta_{ij}:M_i\to M_j\}$ be a direct system of graded $S$-modules where $\theta_{ij}$ are degree-preserving $S$-module homomorphisms. Then
$$\sideset{^*}{_S}\Hom(\varinjlim_iM_i,N)\cong \sideset{^*}{}\varprojlim_i\sideset{^*}{_S}\Hom(M_i,N).$$
\end{pro}
\begin{proof}[Proof]
It suffices to show that
$$\Hom_S(\varinjlim_iM_i,N)_n\cong \varprojlim_i(\Hom_S(M_i,N)_n)$$
for every integer $n$.\par
Let $\theta_i$ denote the natural (degree-preserving) map $M_i\to \varinjlim_iM_i$. Define $$g:\Hom_S(\varinjlim_iM_i,N)_n\to
\varprojlim_i(\Hom_S(M_i,N)_n)$$ as follows. For any
$\alpha\in\Hom_S(\varinjlim_iM_i,N)_n$,
$g(\alpha):=(\cdots,\alpha\circ\theta_i,\cdots)$. Since
$\alpha\circ\theta_j\circ\theta_{ij}=\alpha\circ\theta_i$ and the transition
homomorphisms in the inverse system $\{\Hom_S(M_i,N)_n\}$ are given by composing
with $\theta_{ij}$ which are degree-preserving, $g$ is well-defined.\par
First we prove that $g$ is an injection.  Assume that $\alpha\in \Hom_S(\varinjlim_iM_i,N)_n$ is not 0 and $g(\alpha)=0$. Then there exists a nonzero element $y\in \varinjlim_iM_i$ such that $\alpha(y)\neq 0$. Since $\{M_i\}$ is a direct system, there exist an index $i$ and an element $y_i\in M_i$ such that $\theta_i(y_i)=y$. Then
$\alpha\circ\theta_i(y_i)=\alpha(y)\neq 0$, hence $\alpha\circ\theta_i\neq 0$, so
$g(\alpha)\neq 0$, a contradiction.\par
Next we prove that $g$ is a surjection. Let $(\cdots,\alpha_i,\cdots)$ be an
arbitrary element in $\varprojlim_i(\Hom_S(M_i,N)_n)$, i.e., $\alpha_i$ satisfy
$\alpha_i=\alpha_j\circ\theta_{ij}$ for all $j\geq i$. From the universal property of $\varinjlim_iM_i$, there exists a (unique) homomorphism
$\alpha\in \Hom_S(\varinjlim_iM_i, N)_n$ such that $\alpha\circ\theta_i=\alpha_i$, hence
$g(\alpha)=(\cdots,\alpha_i,\cdots)$.
\end{proof}

We end this section with two basic results on local cohomology that we will use repeatedly.
\begin{itemize}
\item[(GLD)] \textbf{Graded Local Duality.} Let $M$ be a finite graded $R$-module, then there is a functorial isomorphism
$$H^t_{\fm}(M)_l\cong \Hom_k(\Ext^{n+1-t}_R(M,\Omega)_{-l},k)$$
for all integers $t$ and $l$ (cf. \cite[13.4.6]{bs}).\par
If we use $\D(-)$ to denote the graded Matlis dual (for any graded $R$-module $M$, the graded Matlis dual $\D(M)$ of $M$ is defined to be the graded $R$-module with $\D(M)_l=\Hom_k(M_{-l},k)$, then Graded Local Duality can be stated as
\[H^t_{\fm}(M)\cong \D(\Ext^{n+1-t}_R(M,R(-n-1))),\]
for all finitely generated graded $R$-modules $M$ (see \cite[13.3 and 13.4]{bs} for details), and the isomorphism is degree-preserving.

\item[(LCSC)] \textbf{Connection between local cohomology and sheaf cohomology.} Let $M$ be a finite graded $R$-module and let $\widetilde{M}$ be the sheaf on $\mathbb{P}^n$ associated to $M$. Then (cf. \cite[A4.1]{e}) there are a functorial isomorphism
$$H^t_{\fm}(M)\cong \bigoplus_{l\in \mathbb{Z}}H^{t-1}(\mathbb{P}^n,\widetilde{M}(l))\ {\rm when}\ t\geq 2$$
and an exact sequence (functorial in $M$)
$$0\to H^0_{\fm}(M)\to M\to \bigoplus_{l\in \mathbb{Z}}H^0(\mathbb{P}^n,\widetilde{M}(l))\to H^1_{\fm}(M)\to 0$$
where all maps are degree-preserving.
\end{itemize}
\section{The Functors $\E^{s,t}$}
In this section we introduce functors $\E^{s,t}$ from the category of quasi-coherent sheaves on a scheme $Z$ to the category of $\Gamma(Z,\mathcal{O}_Z)$-modules and prove that this functor commutes with finite morphisms (Definition \ref{E-ij} and Corollary \ref{f-finite-gamma}). We start with the following definition.\par

\begin{defn}
\label{omega-z}
If $Z$ is a scheme of finite type over a field $k$ with structure morphism $f:Z\to \Spec(k)$, we define $\dualz$ to be $f^!\mathcal{O}_{\Spec(k)}$.
\end{defn}

\begin{rem}[Uniqueness of dualizing complexes]
\label{unique-dualizing-complex}
{\rm If both $\mathcal{R}^{\bullet}$ and $\mathcal{R}^{\bullet}_1$ are dualizing complexes on a scheme $Y$, then by \cite[Theorem 3.1 in Chapter V]{rd} there exists an invertible sheaf of $\cO_Y$-modules $\mathcal{L}$ such that 
\[\mathcal{R}^{\bullet}\cong\mathcal{R}^{\bullet}_1\otimes\mathcal{L}[n]\]
for some integer $n$, where $[n]$ means a shift to the left by $n$ places (see Section 1 for details).}
\end{rem}

The next lemma shows that, under mild restrictions, the homological shift  $n$ in Remark \ref{unique-dualizing-complex} is 0.
\begin{lem}
\label{dualizing-structure}
Let $X$ be a scheme of finite type over a field $k$ and let $f,g:X\to \Spec(k)$ be two morphisms, then there exists an invertible sheaf $\mathcal{L}$ such that
\[f^!\cO_{\Spec(k)}\cong g^!\cO_{\Spec(k)}\otimes_X\mathcal{L}.\]
\end{lem}
\begin{proof}[Proof]
Without loss of generality, we may and we do assume that $f:X\to \Spec(k)$ is the structure morphism of $X$ over $k$. Hence $f^!\cO_{\Spec(k)}=\dualx$. By Remark \ref{unique-dualizing-complex}, we have $g^!\cO_{\Spec(k)}\cong \dualx \otimes \mathcal{L}[n]$ for some invertible sheaf $\mathcal{L}$ and integer n. We need to show that the shift $n$ is 0. Let $x\in X$ be a closed point such that $k(x)$, the residue field of the local ring $\mathcal{O}_{X,x}$, is a finite extension of $k$ (since $X$ is of finite type over $k$, such a point always exists) and let $i:x\to X$ be the natural embedding. It is clear that the composition $x\xrightarrow{i}X\xrightarrow{f}\Spec(k)$ is the structure morphism $x\xrightarrow{f_x}\Spec(k)$ of $x$, i.e., a finite map $\Spec(k(x))\to \Spec(k)$. Let $U$ be an open neighborhood of $x$ in $X$ such that $\mathcal{L}|_U\cong \cO_X|_U$. Then $i:x\to X$ is the composition of $x\hookrightarrow U\hookrightarrow X$ and hence one has $i^!(\dualx) \otimes \mathcal{L})=i^!\dualx=f^!_x(\mathcal{O}_{\Spec(k)})=k(x)$. It follows that $i^!(g^!\cO_{\Spec(k)})=k(x)[n]$. That is, $n$ is the degree in which the non-zero homology of $i^!(g^!\cO_{\Spec(k)})$ is concentrated.\par
We have a commutative diagram:
$$\xymatrix{x \ar[r]^i \ar[dr]^h & X\ar[d]^g \\
  & \Spec(k)}$$
where $h=g|_{\{x\}}$. The commutativity of the diagram implies that $i^!(g^!\cO_{\Spec(k)})=h^!(\cO_{\Spec(k)})$. But $h$ is a finite map between the spectra of two fields and it is straightforward that $h^!(\cO_{\Spec(k)})=k(x)$, hence the non-zero homology is concentrated in degree zero.
\end{proof}

In particular, if $X,Y$ are $k$-schemes and $g:X\to Y$ is a finite morphism (not necessarily a $k$-morphism), then there exists an invertible sheaf $\mathcal{L}$ on $X$ such that
\[g^!\dualy\cong\dualx\otimes \mathcal{L}.\]

\begin{defn}
\label{E-ij}
For any scheme $Z$ of finite type over a field $k$, we define $$\E^{s,t}(Z;\mathcal{F})\stackrel{{\rm def}}{=}\Ext^s_Z(\underline{\Ext}^t_Z(\mathcal{F},\dualz),\dualz)$$ for any coherent $\mathcal{O}_Z$-sheaf $\mathcal{F}$, where $\dualz$ is as in Definition \ref{omega-z}.
\end{defn}

The main result of this section is the following 
\begin{cor}
\label{f-finite-gamma}
Let $X$ and $Y$ be schemes of finite type over a field $k$. Let $f:X\to Y$ be a finite morphism. For any coherent $\mathcal{O}_X$-sheaf $\mathcal{F}$, there is a natural isomorphism of $\Gamma(Y)$-modules
$$\E^{s,t}(Y;f_*\mathcal{F})\cong f^{\Gamma}_*\E^{s,t}(X;\mathcal{F})$$
\end{cor}
\begin{proof}[Proof]
As a direct consequence from Corollary \ref{duality-sheaf-ext}, one has\footnote{Note that here Corollary \ref{duality-sheaf-ext} is applied twice and one can do this since, by \cite[Proposition 3.3 in Chapter II]{rd}, $\underline{\Ext}^t_X(\mathcal{F},f^!\dualy)$ is coherent.}
$$\E^{s,t}(Y;f_*\mathcal{F})\cong f^{\Gamma}_*\Ext^s_X(\underline{\Ext}^t_X(\mathcal{F},f^!\dualy),f^!\dualy),$$
hence, it suffices to show that
$$\Ext^s_X(\underline{\Ext}^t_X(\mathcal{F},f^!\dualy),f^!\dualy)\cong \E^{s,t}(X;\mathcal{F}).$$
Since $f^!\dualy$ is also a dualizing complex on $X$, thus, according to the remark following Lemma \ref{dualizing-structure}, there exists an invertible sheaf $\mathcal{L}$ on $X$ such that
$$f^!\dualy\cong\dualx\otimes \mathcal{L}.$$
Therefore,
\begin{align}
\Ext^s_X(\underline{\Ext}^t_X(\mathcal{F},f^!\dualy),f^!\dualy) & \cong \Ext^s_X(\underline{\Ext}^t_X(\mathcal{F},\dualy\otimes\mathcal{L}),\dualx\otimes\mathcal{L})\notag\\
&\cong \Ext^s_X(\underline{\Ext}^t_X(\mathcal{F},\dualx)\otimes\mathcal{L},\dualx\otimes\mathcal{L})\tag{*}\\
&\cong \Ext^s_X(\underline{\Ext}^t_X(\mathcal{F},\dualx),\dualx\otimes\mathcal{L}\otimes\mathcal{L}^{-1})\tag{**}\\
&\cong \Ext^s_X(\underline{\Ext}^t_X(\mathcal{F},\dualx),\dualx)=\E^{s,t}(X;\mathcal{F})\notag
\end{align}
where (*) and (**) are direct consequences of Proposition 6.7 in Chapter 3 of \cite{ag} which says that, over any scheme $Z$, one has
$$\Ext^l_Z(\mathcal{F}\otimes \mathcal{L},\mathcal{G})\cong \Ext^l_Z(\mathcal{F},\mathcal{L}^{-1}\otimes \mathcal{G})$$ and
$$\underline{\Ext}^l_Z(\mathcal{F}\otimes \mathcal{L},\mathcal{G})\cong \underline{\Ext}^l_Z(\mathcal{F},\mathcal{L}^{-1}\otimes \mathcal{G})\cong \underline{\Ext}^l_Z(\mathcal{F}, \mathcal{G})\otimes\mathcal{L}^{-1}$$
where $\mathcal{F}$ and $\mathcal{G}$ are sheaves of $\mathcal{O}_Z$-modules, $\mathcal{L}$ is an invertible sheaf of $\mathcal{O}_Z$-modules and $\mathcal{L}^{-1}$ is the dual of $\mathcal{L}$.
\end{proof}

\section{The $k$-vector spaces $\mathcal{M}^{i,j}_0$}
In this section we will prove Theorem \ref{indp-M-i-j} and some related results. Let $X$ be a projective scheme over a field $k$ and let $\eta:X\hookrightarrow \mathbb{P}^n_k$ be an embedding. Let $R$ denote $k[x_0,\dots,x_n]$, the coordinate ring of $\mathbb{P}^n_k$, and let $I\subset R$ denote the defining ideal of $X$ under this embedding. It is clear that $\widetilde{\Omega}=\omega_{\mathbb{P}^n}$, where $\widetilde{\Omega}$ denotes the sheaf associated to $\Omega$ on $\mathbb{P}^n$. \par

\begin{lem}
\label{lemma-degree-0-ext-sheaf}
Let $M$ be a finitely generated $\mathbb{Z}$-graded $R$-module. Then there is a functorial $k$-linear map for any $i\geq 1$
\begin{equation}
\label{degree-0-ext-sheaf-ext}
\Ext^{n+1-i}_R(M,\Omega)_0\to \Ext^{1-i}_{\mathbb{P}^n}(\widetilde{M},\dualp)
\end{equation}
which is an isomorphism when $i\geq 2$. The functoriality here means that, if there is a degree-preserving $R$-module homomorphism $M\to N$, one has a commutative diagram
\[\xymatrix{
\Ext^{n+1-i}_R(M,\Omega)_0\ar[d] & \Ext^{n+1-i}_R(N,\Omega)_0\ar[l]\ar[d]\\
\Ext^{1-i}_{\mathbb{P}^n}(\widetilde{M},\dualp) & \Ext^{1-i}_{\mathbb{P}^n}(\widetilde{N},\dualp)\ar[l]\ .
}
\]
\end{lem}
\begin{proof}[Proof]
The map (\ref{degree-0-ext-sheaf-ext}) is given by 
\begin{align}
&\quad \Ext^{n+1-i}_R(M,\Omega)_0\notag\\
&\xrightarrow{\sim}\Hom_k(H^i_{\fm}(M)_0,k)\ {\rm by\ (GLD)}\notag\\
&\to \Hom_k(H^{i-1}(\mathbb{P}^n,\widetilde{M}),k)\tag{\dag}\\
&\xrightarrow{\sim}\Ext^{n+1-i}_{\mathbb{P}^n}(\widetilde{M},\omega_{\mathbb{P}^n})\ {\rm by\ Serre\ Duality}\notag\\
&=\Ext^{1-i}_{\mathbb{P}^n}(\widetilde{M},\dualp)\tag{\ddag} 
\end{align}
where the map ($\dag$) follows from (LCSC) and is an isomorphism when $i\geq 2$; and ($\ddag$) holds because $\dualp=\omega_{\mathbb{P}^n}[n]$ (a straightforward consequence from the definition on page 145 in \cite{rd}). Recall that $\omega_{\mathbb{P}^n}[n]$ means the complex with $\omega_{\mathbb{P}^n}$ in homological degree $-n$ and zero elsewhere. Since every map involved is functorial, so is (\ref{degree-0-ext-sheaf-ext}). This completes the proof.
\end{proof}

For any $k$-vector space $V$ we will use $\D(V)$ to denote its dual space $\Hom_k(V,k)$ assuming no confusion will arise.\par
  
The following is the main result of this section.
\begin{pro}
\label{M-ij}
Let $L$ be a module-finite graded $R/I$-algebra such that in degree 0 the natural map $k=(R/I)_0\to L_0$ is a bijection. Let $\widetilde{L}$ be the sheaf on $X$ associated to $L$. Then
\begin{enumerate}
\item when $i\geq 2$, there is an isomorphism of $k$-vector spaces
\begin{equation}
\label{L-ij-2}
\cE^{i,j}(L)_0\cong \E^{1-i,1-j}(X;\widetilde{L})
\end{equation}
which is functorial in $L$
\item when $i=0,1$, one has the following exact sequence of $k$-vector spaces 
\begin{equation}
\label{L-ij-0-1}
0\to \cE^{1,j}(L)_0\to \E^{0,1-j}(X;\widetilde{L})\xrightarrow{\delta} H^j_{\fm}(L)_0\to \cE^{0,j}(L)_0 \to 0,
\end{equation}
such that $H^j_{\fm}(L)_0$ and $\delta$ only depend on the sheaf of $\cO_X$-modules $\widetilde{L}$ (it is evident that $\E^{0,1-j}(X;\widetilde{L})$ only depends on $\widetilde{L}$). This sequence is functorial in $L$.
\end{enumerate}
\end{pro}

Before we proceed to the proof, we explain why $L$ is assumed to be a finite $R/I$-module. The module-finiteness of $L$ over $R/I$ ensures that, during the course of the proof of Proposition \ref{M-ij} all the $R/I$-modules to which (GLD) will be applied are finite graded $R$-modules; this is where the module-finiteness of $L$ is used.

\begin{proof}[Proof]
Let $g:\mathbb{P}^n_k\to \Spec(k)$ denote the structure morphism and let $\dualp$ denote $g^!\mathcal{O}_{\Spec(k)}$. Then $\dualp=\omega_{\mathbb{P}^n_k}[n]$. When $i\geq 2$, one has \begin{align}
\cE^{i,j}(L)_0 &= \Ext^{n+1-i}_R(\Ext^{n+1-j}_R(L,\Omega),\Omega)_0\\
&\cong \Ext^{1-i}_{\mathbb{P}^n}(\underline{\Ext}^{1-j}_{\mathbb{P}^n}(\eta_*\widetilde{L},\dualp),\dualp)  \tag{i}\\
&\cong \eta^{\Gamma}_*\E^{1-i,1-j}(X;\widetilde{L})\tag{ii}\\
&=\E^{1-i,1-j}(X;\widetilde{L})\tag{iii}
\end{align}
(i) follows from applying (\ref{degree-0-ext-sheaf-ext}) to $M=\Ext^{n+1-j}_R(L,\Omega)$ in the case when $i\geq 2$ and Proposition \ref{sheaf-ext}. It is here that we use $i\geq 2$.\\
(ii) is a consequence of Corollary \ref{f-finite-gamma}.\\
(iii) is true since $\eta^{\Gamma}_*\E^{1-i,1-j}(X;\widetilde{L})$ and $\E^{1-i,1-j}(X;\widetilde{L})$ are identical as $k$-vector spaces.\par

Since isomorphisms (i) and (ii) are functorial in $L$, the isomorphism \[\mathcal{L}^{i,j}_0\cong \E^{1-i,1-j}(X;\widetilde{L})\] itself is functorial in $L$. This completes the proof of part (1) in Proposition \ref{M-ij}.\par 

Consider the module $T=\Ext^{n+1-j}_R(L,\Omega)$ with $\widetilde{T}=\underline{\Ext}^{1-j}_{\mathbb{P}^n}(\eta_*\widetilde{L},\dualp)$ (by Proposition 2.3). From (LCSC), we get the following exact sequence 
\begin{equation}
\label{g-s-0-1}
0\to H^0_{\mathfrak{m}}(T)_0\to T_0\xrightarrow{\alpha} H^0(\mathbb{P}^n,\widetilde{T})\to H^1_{\mathfrak{m}}(T)_0\to 0
\end{equation}
which induces an exact sequence
$$0\to \D(H^1_{\mathfrak{m}}(T)_0)\xrightarrow{h_1} \D(H^0(\mathbb{P}^n,\widetilde{T}))\xrightarrow{h_2} \D(T_0) \xrightarrow{h_3} \D(H^0_{\mathfrak{m}}(T)_0)\to 0,$$
where $\D(-)$ denotes the dual space $\Hom_k(-,k)$.\par
Now consider the following commutative diagram (where $\Omega$ denotes $R(-n-1)$)
\begin{equation}
\label{construction-delta}
\xymatrix{
0 \ar[d] & 0 \ar[d]\\
 \D(H^1_{\mathfrak{m}}(T)_0) \ar[r]^{\sim}_{q_1} \ar[d]^{h_1} &  \Ext^n_R(T,\Omega)_0 \ar[d]^{q_2\circ h_1\circ q^{-1}_1}\\
\D(H^0(\mathbb{P}^n,\widetilde{T})) \ar[r]^{\sim}_{q_2} \ar[d]^{h_2} &  \Ext^0_{\mathbb{P}^n}(\widetilde{T},\dualp) \ar[d]^{q_3\circ h_2\circ q^{-1}_2}\\
 \D(T_0) \ar[r]^{\sim}_{q_3} \ar[d]^{h_3} & H^j_{\mathfrak{m}}(L)_0 \ar[d]^{q_4\circ h_3\circ q^{-1}_3}\\
 \D(H^0_{\mathfrak{m}}(T)_0) \ar[r]^{\sim}_{q_4} \ar[d] & \Ext^{n+1}_R(T,\Omega)_0 \ar[d]\\
0& 0 }
\end{equation}
where $q_1$, $q_3$ and $q_4$ are consequences of (GLD) and $q_2$ follows from the Serre Duality (Theorem 7.1 in Chapter 3 of \cite{ag}) plus the fact that $\dualp=\omega_{\mathbb{P}^n}[n]$; the maps in the right column are induced by the maps in the left column.\par
Since $\Ext^n_R(T,\Omega)_0=\cE^{1,j}(L)_0$, while $\Ext^{n+1}_R(T,\Omega)_0=\cE^{0,j}(L)_0$, and 
\begin{equation}
\label{E-0-1-j}
\Ext^0_{\mathbb{P}^n}(\widetilde{T},\dualp)=\E^{0,1-j}(\mathbb{P}^n;\eta_*\widetilde{L})\stackrel{{\rm Cor.\ \ref{f-finite-gamma}}}{\cong} \eta^{\Gamma}_*\E^{0,1-j}(X;\widetilde{L})=\E^{0,1-j}(X;\widetilde{L})
\end{equation}
the right column of the commutative diagram (\ref{construction-delta}) induces the desired exact sequence
$$0\to \cE^{1,j}(L)_0\to \E^{0,1-j}(X;\widetilde{L})\xrightarrow{\delta} H^j_{\mathfrak{m}}(L)_0\to \cE^{0,j}(L)_0 \to 0.$$
Since every step involved in constructing this sequence is functorial in $L$, the sequence itself is functorial in $L$.\par
It remains to show that $H^j_{\fm}(L)_0$ and $\delta$ only depend on $\widetilde{L}$.\par
First, we will prove that $H^j_{\fm}(L)_0$ only depends on $\widetilde{L}$. If $\widetilde{L}=0$, then $H^0_{\fm}(L)=L$ and 
$H^j_{\fm}(L)=0$ for $j>0$, i.e. $H^0_{\fm}(L)_0=k$ and $H^j_{\fm}(L)_0=0$ for $j>0$. Thus the condition $\widetilde{L}=0$ determines $H^j_{\fm}(L)_0$. It remains to consider the case $\widetilde{L}\neq 0$. In this case, we consider two cases: $j\geq 2$ and $j\leq 1$.\par
When $j\geq 2$, we have
$$H^j_{\fm}(L)_0\cong H^{j-1}(\mathbb{P}^n_k,\eta_*\widetilde{L})=H^{j-1}(X,\widetilde{L})$$
by (LCSC) and the fact that $\eta$ is a finite morphism. This certainly only depends on $\widetilde{L}$.\\
When $j\leq 1$, we consider the exact sequence from (LCSC)

\begin{equation}
\label{T-0-1}
{\tiny
\xymatrix{
0\ar[r]  & H^0_{\fm}(L)_0\ar[r] & L_0\ar[r] \ar@{=}[d]&H^0(\mathbb{P}^n_k,\eta_*\widetilde{L})\ar[r] \ar@{=}[d]&  H^1_{\fm}(L)_0\ar[r] & 0\\
 & & k & H^0(X,\widetilde{L})& & .}}
\end{equation}

The condition $\widetilde{L}\neq 0$ implies $1\notin H^0_{\fm}(L)$. Since $H^0_{\fm}(L)_0\subset L_0=k$ and $L_0$ is the $k$-linear 
span of $1$ and $1\notin H^0_{\fm}(L)$, we see that $H^0_{\fm}(L)_0=0$, therefore,
$H^1_{\fm}(L)_0\cong H^0(X,\widetilde{L})/k$ which only depends on $\widetilde{L}$. This shows that $H^j_{\fm}(L)_0$ only depends on 
$\widetilde{L}$.\par
Finally we prove that $\delta$ only depends on $\widetilde{L}$. If $\widetilde{L}=0$, then 
$\E^{0,1-j}(X;\widetilde{L})=0$, hence $\delta=0$, i.e. $\delta$ is determined by $\widetilde{L}$. It 
remains to consider the case $\widetilde{L}\neq 0$.\par
From the construction of $\delta$ (the diagram (\ref{construction-delta}) and the isomorphism (\ref{E-0-1-j})), it suffices to show that the map $h_2:\D(H^0(\mathbb{P}^n,\widetilde{T}))\to \D(T_0)$ in the diagram (\ref{construction-delta}) depends only on $\widetilde{L}$ and hence it is enough to show that the dual map $\alpha:T_0\to H^0(\mathbb{P}^n,\widetilde{T})$ in the exact sequence (\ref{g-s-0-1}) depends only on $\widetilde{L}$.\par
A proof of this fact is unexpectedly long. For a graded ring $S$ and a graded $S$-module $M$, an element of $M_0$
automatically produces a global section of $\widetilde{M}$ on $\Proj(S)$. The map $\alpha$ is to send each element of $T_0$ to its
associated global section of $\widetilde{T}$ on $X$.\par 
Since $\widetilde{L}\neq 0$ implies $H^0_{\fm}(L)_0=0$, i.e., $T_0\cong \Hom_k(H^0_{\fm}(L)_0,k))=0$ when $j=0$, we may assume that $j\geq 1$.\par

One always has
\begin{align}
\widetilde{T} &= \eta^*\underline{\Ext}^{n+1-j}_{\mathbb{P}^n}(\eta_*\widetilde{L},\omega_{\mathbb{P}^n})\notag\\
&\cong \eta^*\underline{\Ext}^{1-j}_{\mathbb{P}^n}(\eta_*\widetilde{L},\dualp)\notag\\
&\cong \eta^*\eta_*\underline{\Ext}^{1-j}_X(\widetilde{L},\eta^!\dualp)\notag\\
&=\underline{\Ext}^{1-j}_X(\widetilde{L},\dualx)\notag
\end{align}
From (LCSC), the natural map
$H^{j-1}(X,\widetilde{L})=H^{j-1}(\mathbb{P}^n_k,\eta_*\widetilde{L})\to
H^{j}_{\fm}(L)_0$, an isomorphism when $j\geq 2$ and a
surjection when $j=1$, induces
$$\xi:\Hom_k(H^{j}_{\fm}(L)_0,k)\to \Hom_k(H^{j-1}(\mathbb{P}^n,\eta_*\widetilde{L}),k)$$
which is an isomorphism when $j\geq 2$ and an injection when $j=1$.
Hence, one always has a map 
\begin{equation}
\label{theta}
\theta:T_0\to \Ext^{1-j}_X(\widetilde{L},\dualx)
\end{equation}
as follows
\begin{align}
T_0 &= \Ext^{n+1-j}_R(L,\Omega)_0\notag\\
&\cong \Hom_k(H^j_{\fm}(L)_0,k)\notag\\
&\xrightarrow{\xi}
 \Hom_k(H^{j-1}(\mathbb{P}^n,\eta_*\widetilde{L}),k)(=\Hom_k(H^{j-1}(X,\widetilde{L}),k))\notag\\
&\cong \Ext^{n+1-j}_{\mathbb{P}^n}(\eta_*\widetilde{L},\omega_{\mathbb{P}^n})\tag{\dag}\\
&=\Ext^{1-j}_{\mathbb{P}^n}(\eta_*\widetilde{L},\dualp)\notag\\
&\cong \Ext^{1-j}_X(\widetilde{L},\dualx)\tag{\ddag}
\end{align}
where $(\dag)$ is the inverse of the isomorphism in Corollary 3.4(c) in Chapter VII in \cite{rd} applied to the structure morphism $g:\mathbb{P}^n\to k$ and $(\ddag)$ is the inverse of the isomorphism that follows from Corollary \ref{duality-sheaf-ext} applied to the embedding $\eta:X\to \mathbb{P}^n$ with $\mathcal{F}=\widetilde{L}$ and $\mathcal{G}^{\bullet}=\dualx$ (note that $\eta^{\Gamma}_*\Ext^{1-j}_X(\widetilde{L},\dualx)$ is, as a $k$-vector space, identical as $\Ext^{1-j}_X(\widetilde{L},\dualx)$ ).\par
Since $\xi$ is an isomorphism for $j\geq 2$ and injective for $j=1$,
so is $\theta$. When $j=1$, we consider the exact sequence (\ref{T-0-1}) again
$$0\to k\to H^0(X,\widetilde{L})\to H^1_{\mathfrak{m}}(L)_0\to 0.$$
The map $k\to H^0(X,\widetilde{L})$ sends each element $c\in k$ to a
string $(\dots,c_t,\dots)$ (which is a section of $\widetilde{L}$ on $X$) with $c_t$ at each connected component $X_t$ of
$X$ where $c_t=0$ if $\widetilde{L}|_{X_t}=0$ and $c_t=c$ if $\widetilde{L}|_{X_t}\neq 0$. Hence $k\to H^0(X,\widetilde{L})$ is injective and depends only on $\widetilde{L}$. But the map $\xi$ is nothing but the dual of $H^0(X,\widetilde{L})\to H^1_{\mathfrak{m}}(L)_0$. Hence the image of the map $\xi$ is the kernel of the dual of $k\to H^0(X,\widetilde{L})$ and therefore depends only on $\widetilde{L}$. The composition of $(\dag)$ and $(\ddag)$ in the definition of $\theta$ depends only on $\widetilde{L}$ because it is the inverse of the isomorphism in Corollary 3.4(c) in Chpater VII in \cite{rd} \footnote{This holds because it is clear that the structure morphism $X\to \Spec(k)$ satisfies $f=g\circ\eta$ and hence$f^!=(g\circ\eta)^!=\eta^!g^!$ by \cite[Theorem 8.7 in Chapter III]{rd}.}
with $f$ being the structure morphism $X\to \Spec(k)$, $F^{\bullet}=\widetilde{L}$, and $G^{\bullet} =\cO_{\Spec(k)}$. When $j\geq 2$, we have seen that the map $\theta$ is an isomorphism. Thus the map $\theta$ maps $T_0$ bijectively onto a subspace of $\Ext^{1-j}_X(\widetilde{L},\dualx)$ that only depends on $\widetilde{L}$ (the codimension of this subspace in $\Ext^{1-j}_X(\widetilde{L},\dualx)$ is 0 or 1 depending on whether or not $j\geq 2$).\par

Let $\mathcal{I}^{\bullet}$ be an injective
resolution of $\dualx$, then $\Ext^{t}_X(\widetilde{L},\dualx)$ and $\underline{\Ext}^{t}_X(\widetilde{L},\dualx)$ may be calculated by the
following two complexes for all integers $t$
$$\cdots\to \Hom_X(\widetilde{L},\mathcal{I}^{t-1})\xrightarrow{d_{t-1}}\Hom_X(\widetilde{L},\mathcal{I}^t)\xrightarrow{d_t}\Hom_X(\widetilde{L},\mathcal{I}^{t+1})\to\cdots$$
and
$$\cdots\to \underline{\Hom}_X(\widetilde{L},\mathcal{I}^{t-1})\xrightarrow{\delta_{t-1}}\underline{\Hom}_X(\widetilde{L},\mathcal{I}^t)\xrightarrow{\delta_t}\underline{\Hom}_X(\widetilde{L},\mathcal{I}^{t+1})\to\cdots$$
$\underline{\Ext}^t_X(\widetilde{L},\dualx)$ is the sheaf associated to
the presheaf $U\mapsto \ker(\delta_t)(U)/\im(\delta_{t-1})(U)$, hence there is a map $$\Gamma(X,\ker(\delta_t)/\im(\delta_{t-1}))\to
\Gamma(X,\underline{\Ext}^{t}_X(\widetilde{L},\dualx)).$$ It is clear that
there is a map $$\Ext^t_X(\widetilde{L},\dualx)=\ker(d_t)/\im(d_{t-1})\to
\Gamma(X,\ker(\delta_t)/\im(\delta_{t-1})).$$ The composition of
these two maps produces a natural map
$$\Ext^t_X(\widetilde{L},\dualx)\xrightarrow{\gamma} \Gamma(X,\underline{\Ext}^{t}_X(\widetilde{L},\dualx))$$
for all integers $t$, which clearly depends only on $\widetilde{L}$.\par

Consider the diagram 
\begin{equation}
\label{diagram-sharp}
\xymatrix{
T_0 \ar[r]^{\alpha} \ar[d]^{\theta} & \Gamma(X,\widetilde{T})\ar[d]^{\sim}_u\\
\Ext^{1-j}_X(\widetilde{L},\dualx)\ar[r]^{\gamma} & \Gamma(X,\underline{\Ext}^{1-j}_X(\widetilde{L},\dualx))
}
\end{equation}

The right vertical map is an isomorphism because $\widetilde{T}=\underline{\Ext}^{1-j}_X(\widetilde{L}, \dualx)$, the map $\gamma$ only depends on $\widetilde{L}$ and $\theta$ sends $T_0$ bijectively onto a subspace of $\Ext^{1-j}_X(\widetilde{L},\dualx)$ that only depends on $\widetilde{L}$. Therefore the commutativity of the diagram would show that the map $\alpha$ only depends on $\widetilde{L}$. Thus to complete the proof of Proposition \ref{M-ij} it remains to show that diagram (\ref{diagram-sharp}) is commutative. This is shown below.\par

Let $M$ be a graded $R$-module, then the map $M\to
\oplus_{l\in\mathbb{Z}}\Gamma(\mathbb{P}^n,\tilde{M}(l))$ (by (LCSC)) is a degree-preserving functorial
$R$-module homomorphism.\par

We construct a graded $R$-module $Q$ as follows,
$$Q_l:=\Ext^{n+1-j}_{\mathbb{P}^n}(\eta_*\widetilde{L}(-l),\omega_{\mathbb{P}^n}).$$
When $j\geq 1$, we have a natural map (by applying $\Hom_k(-,k)$ to (LCSC))
$$T_l\cong \Hom_k(H^j_{\fm}(L)_{-l},k)\xrightarrow{\chi_l} \Hom_k(H^{j-1}(\mathbb{P}^n,\eta_*\widetilde{L}(-l)),k)=Q_l$$
for all integers $l$. Therefore, we have a degree-preserving map $\chi:T\to Q$ of graded
$R$-modules and hence a commutative diagram
\begin{equation}
\label{T-Q}
\xymatrix{
T\ar[r]\ar[d]^{\chi} & \oplus_{l\in\mathbb{Z}}\Gamma(\mathbb{P}^n,\tilde{T}(l))\ar[d] \\
Q\ar[r] &\oplus_{l\in\mathbb{Z}}\Gamma(\mathbb{P}^n,\tilde{Q}(l))
}
\end{equation}
which induces a commutative diagram of sheaves on $\mathbb{P}^n$
\begin{equation}
\label{T-Q-sheaf-diagram}
\xymatrix{
\tilde{T}\ar[r]\ar[d] & (\oplus_{l\in\mathbb{Z}}\Gamma(\mathbb{P}^n,\tilde{T}(l)))^{\widetilde{ }}\ar[d] \\
\tilde{Q}\ar[r] & (\oplus_{l\in\mathbb{Z}}\Gamma(\mathbb{P}^n,\tilde{Q}(l)))^{\widetilde{ }}}
\end{equation}

We claim that the map $\chi:T\to Q$ induces an isomorphism $\tilde{T}\to \tilde{Q}$ and we reason as follows. According to (LCSC) we have the following isomorphism and exact sequence
$$\oplus_lH^{j-1}(\mathbb{P}^n,\widetilde{L}(l))\xrightarrow{\sim}H^j_{\fm}(L),\ {\rm when}\ j\geq 2$$
$$0\to H^0_{\fm}(L)\to L\to \oplus_lH^0(\mathbb{P}^n,\eta_*\widetilde{L}(l))\to H^1_{\fm}(L)\to 0$$

which in turn induce (after applying the graded Matlis dual (\cite[13.4.5]{bs}))
\begin{enumerate}
\item[(i)] an isomorphism $T\xrightarrow{\chi}Q$ when $j\geq 2$; 
\item[(ii)] an exact sequence
$$0\to T\xrightarrow{\chi}Q\to \oplus_l\Hom_k((L)_{-l},k)\ {\rm when}\ j=1$$
\end{enumerate}
From (i), we can see that when $j\geq 2$ the map $\chi:T\to Q$ is actually an isomorphism; hence induces an isomorphism $\tilde{T}\to \tilde{Q}$. When $j=1$, (ii) induces an exact sequence of the associated sheaves on $\mathbb{P}^n$
$$0\to \tilde{T}\to \tilde{Q}\to (\oplus_l\Hom_k(L_{-l},k))^{\widetilde{ }}$$
Set $N=\oplus_l\Hom_k(L_{-l},k)$. Then $N$ is the graded Matlis dual of $L$ (cf. \cite[13.4.5]{bs}) which is clearly only supported at $\fm$; and hence the associated sheaf $\tilde{N}$ on $\mathbb{P}^n$ is 0. This proves that $\chi$ induces an isomorphism $\tilde{T}\cong \tilde{Q}$ when $j=1$ and hence the proof of our claim.\par
From  \cite[Proposition 5.15 in ChapterII]{ag}, we know that the two horizontal maps in the diagram (\ref{T-Q-sheaf-diagram}) are isomorphisms and hence all maps in (\ref{T-Q-sheaf-diagram}) are isomorphisms. Therefore the diagram (\ref{T-Q}) induces the following commutative diagram,
$${\tiny
\xymatrix{
 & T_0\ar[r]\ar[d] & \Gamma(\mathbb{P}^n,\tilde{T})\ar[d] & \\
\Ext^{n+1-j}(\eta_*\widetilde{L},\omega_{\mathbb{P}^n})\ar@{=}[r] &
Q_0\ar[r] & \Gamma(\mathbb{P}^n,\tilde{Q})\ar@{=}[r] &
\Gamma(\mathbb{P}^n,\tilde{T}=\underline{\Ext}^{n+1-j}_{\mathbb{P}^n}(\eta_*\widetilde{L},\omega_{\mathbb{P}^n}))
}}$$ 
where $\tilde{T}=\underline{\Ext}^{n+1-j}_{\mathbb{P}^n}(\eta_*\widetilde{L},\omega_{\mathbb{P}^n})$ holds because of Proposition 2.3.\par
Diagram (\ref{diagram-sharp}) can be written as
$$\xymatrix{
T_0 \ar[r]\ar[d] & \Gamma(\mathbb{P}^n,\tilde{T})\ar@{=}[d] \\
\Ext^{n+1-j}_{\mathbb{P}^n}(\eta_*\widetilde{L},\omega_{\mathbb{P}^n})\ar[r]\ar@{=}[d] & \Gamma(\mathbb{P}^n,\underline{\Ext}^{n+1-j}_{\mathbb{P}^n}(\eta_*\widetilde{L},\omega_{\mathbb{P}^n}))\ar@{=}[d] \\
\Ext^{1-j}_{\mathbb{P}^n}(\eta_*\widetilde{L},\dualp) \ar[r]\ar[d]^{\sim} & \Gamma(X,\underline{\Ext}^{1-j}_{\mathbb{P}^n}(\eta_*\widetilde{L},\dualp))\ar[d]^{\sim}\\
\Ext^{1-j}_X(\widetilde{L},\dualx)\ar[r] &
\Gamma(X,\underline{\Ext}^{1-j}_X(\widetilde{L},\dualx)) 
}$$
We already proved the commutativity of the square at the top; the commutativity of the
remaining two squares is clear. Hence diagram (\ref{diagram-sharp}) is
commutative. This completes the proof that $\alpha$ only depends on $\widetilde{L}$ and the proof of Proposition \ref{M-ij}.
\end{proof}

\begin{cor}
\label{independence-L-ij-0}
Let $L$ be be an $R/I$-algebra as in Proposition \ref{M-ij}. Then $\cE^{i,j}(L)_0$ only depends on the sheaf of $\cO_X$-modules $\widetilde{L}$.
\end{cor}
\begin{proof}[Proof]
When $i\geq 2$, by (\ref{L-ij-2}) $\cE^{i,j}(L)_0$ is isomorphic to $\E^{1-i,1-j}(X;\widetilde{L})$ which clearly only depends on $\widetilde{L}$. When $i\leq 1$, since the two middle terms and the map between them in the exact sequence (\ref{L-ij-0-1}) depend only on $\widetilde{L}$, the two terms at the ends, $\cE^{0,j}(L)_0$ and $\cE^{1,j}(L)_0$, only depend on $\widetilde{L}$.  
\end{proof}

It should be pointed out here that $\widetilde{L}=0$ does not imply $\cE^{i,j}(L)_0=0$ for all $i,j$. It implies $\cE^{i,j}(L)_0=0$ if $i\neq 0$ or $j\neq 0$, but $\cE^{0,0}(L)_0=k$. This still means that $\cE^{i,j}(L)_0$ is determined by $\widetilde{L}$.\par

So far we have proved Theorem \ref{indp-M-i-j} which is a special case of Corollary \ref{independence-L-ij-0} with $L=R/I$. The rest of this section is devoted to proving that various maps are independent of the embedding, which is necessary for Section 7.\par 

Let $L$ be an $R/I$-algebra as in Proposition \ref{M-ij}, then the natural map $R/I\to L$ induces the map $\pi:\cE^{i,j}(R/I)_0\to \cE^{i,j}(L)_0$. It follows from Proposition \ref{M-ij} that $\cE^{i,j}(R/I)_0$ and $\cE^{i,j}(L)_0$ only depend on $\cO_X$ and $\widetilde{L}$ respectively.

\begin{lem}
\label{local-cohomology-independent}
Let $L$ be as in Proposition \ref{M-ij}. The natural map $H^j_{\fm}(R/I)_0\to H^j_{\fm}(L)_0$ only depends on the sheaf of $\cO_X$-modules $\widetilde{L}$.
\end{lem}
\begin{proof}[Proof]
We will consider two cases: $j\geq 2$ and $j\leq 1$.\par
When $j\geq 2$, the following commutative diagram (from (LCSC)) proves that $H^j_{\fm}(R/I)_0\to H^j_{\fm}(L)_0$ only depends on $\widetilde{L}$
$$\xymatrix{
H^j_{\fm}(R/I)_0 \ar[r] \ar[d]^{p_1}_{\sim} & H^j_{\fm}(L)_0 \ar[d]^{q_1}_{\sim}\\
H^{j-1}(\mathbb{P}^n,\eta_*\cO_X) \ar[r] \ar[d]^{p_2}_{\sim} & H^{j-1}(\mathbb{P}^n,\eta_*\widetilde{L}) \ar[d]^{q_2}_{\sim}\\
H^{j-1}(X,\cO_X) \ar[r] & H^{j-1}(X,\widetilde{L})
}$$
where $p_1$ and $q_1$ are from (LCSC) since $j\geq 2$, this is where the assumption $j\geq 2$ is used; $p_2$ and $q_2$ are isomorphisms since $\eta$ is a finite morphism; the second and the third horizontal maps are induced by $\cO_X\to \widetilde{L}$.\par
If $\widetilde{L}=0$, then $H^0_{\fm}(L)_0=L_0=k$ and $H^j_{\fm}(L)=0$ for $j\geq 1$. Hence in this case $H^j_{\fm}(R/I)_0\to H^j_{\fm}(L)_0$ is either the bijection $k\to k$ when $j=0$ or the 0-map when $j\geq 1$. Therefore, $H^j_{\fm}(R/I)_0\to H^j_{\fm}(L)_0$ is determined by $\widetilde{L}$.\par
Assume that $\widetilde{L}\neq 0$. When $j=0$, we have seen from the proof of Proposition \ref{M-ij} that $H^0_{\fm}(R/I)_0=H^0_{\fm}(L)_0=0$, hence there is nothing to prove.\par
Assume that $j=1$. We consider the following commutative diagram (which follows from (LCSC))
{\tiny
$$\xymatrix{
H^0_{\fm}(R/I)_0=0\ar[r] & (R/I)_0=k\ar[r]\ar[d]^{w_1}_{=} & H^0(X,\cO_X)\ar[r]\ar[d]^{w_2} & H^1_{\fm}(R/I)_0\ar[d]_{w_3}\ar[r] & 0\\
H^0_{\fm}(L)_0=0\ar[r] & L_0=k\ar[r] & H^0(X,\widetilde{L})\ar[r] & H^1_{\fm}(L)_0\ar[r] & 0}$$}
where $w_2$ is induced by $\cO_X\to \widetilde{L}$. Therefore, $w_3$ only depends on $\widetilde{L}$. This completes the proof of our lemma.
\end{proof}

\begin{cor}
\label{M-L-0-independence}
The map $\pi:\cE^{i,j}(R/I)_0\to \cE^{i,j}(L)_0$ only depends on $\cO_X$ and the sheaf of $\cO_X$-modules $\widetilde{L}$.
\end{cor}
\begin{proof}[Proof]
When $i\geq 2$, since the isomorphism (\ref{L-ij-2}) is functorial in $L$, we have a commutative diagram
$$\xymatrix{
\cE^{i,j}(R/I)_0 \ar[rr]^{\sim}\ar[d]^{\pi}  & &\E^{1-i,1-j}(X;\cO_X)\ar[d]^{\pi'}\\
\cE^{i,j}(L)_0 \ar[rr]^{\sim} & &\E^{1-i,1-j}(X;\widetilde{L})
}$$
where $\pi'$ is induced by the sheaf morphism $\cO_X\to \widetilde{L}$. It is clear from this commutative diagram that $\pi:\cE^{i,j}(R/I)_0\to \cE^{i,j}(L)_0$ only depends on $\cO_X$ and $\widetilde{L}$ when $i\geq 2$.\par
When $i=0,1$, the functorial exact sequence (\ref{L-ij-0-1}) induces a commutative diagram
{\tiny
$$\xymatrix{
0\ar[r] & \cE^{1,j}(R/I)_0\ar[r]\ar[d] & \E^{0,1-j}(X;\cO_X)\ar[r]\ar[d] & H^j_{\fm}(R/I)_0\ar[r]\ar[d] & \cE^{0,j}(R/I)_0\ar[r]\ar[d] &  0\\
0\ar[r] & \cE^{1,j}(L)_0\ar[r] & \E^{0,1-j}(X;\widetilde{L})\ar[r] & H^j_{\fm}(L)_0\ar[r] & \cE^{0,j}(L)_0\ar[r] &  0
}$$}
where all vertical maps in the middle are induced by $R/I\to L$ (and $\cO_X\to \widetilde{L}$). From part (2) of Proposition \ref{M-ij}, one can see that all entries in this diagram and all horizontal maps depend only on $\cO_X$ and $\widetilde{L}$; it is clear that the second vertical map depends only on $\cO_X$ and $\widetilde{L}$; according to Lemma \ref{local-cohomology-independent} the third vertical map depends only on $\widetilde{L}$; hence the vertical maps at the ends depend only on $\cO_X$ and $\widetilde{L}$ as well. This finishes the proof of the corollary.
\end{proof}

\section{A graded theory of the Frobenius endomorphism}
Throughout Sections 5, 6 and 7, $k$ denotes a perfect field of prime characteristic $p > 0$ and $X$ denotes a projective scheme over $k$. Let $R = k[x_0,\dots,x_n]$ and let $I$ be the defining ideal of $X$ under $\eta:X\hookrightarrow \mathbb{P}^n_k$.\par

In this section we develop a graded version of the theory of the Frobenius endomorphism; this serves as a preparation for the next section where we will construct a graded $p$-linear structure on $\cE^{i,j}(R/I)$. First, we explain what a (graded) $p$-linear structure is.

\begin{defn}
For any $k$-algebra $B$ and a $B$-module $M$, a $p$-linear structure $\varphi$ on $M$ is an endomorphism on the underlying additive group that satisfies $\varphi(bm)=b^p\varphi(m)$ for all $b\in B$ and $m\in M$.\par
When both $B$ and $M$ are graded, a $p$-linear structure $\varphi$ on $M$ is called a graded $p$-linear structure if $\deg(\varphi(m))=p\deg(m)$ for all homogeneous elements $m\in M$.
\end{defn}

The Frobenius endomorphism $F_R:R\to R$ (or simply $F$) is defined by $r\mapsto r^p$. We will not distinguish $F:R\to R$ and the corresponding endomorphism on $\Spec(R)$. The Frobenius pullback $F^*: R{\rm -mod}\to R{\rm -mod}$ is defined as follows.
Let $R^{(1)}$ be the additive group of $R$ regarded as an $R$-bimodule with usual left $R$-action and with right $R$-action given by $r'\cdot r=r^pr'$ for all $r\in R,r'\in R^{(1)}$. The Frobenius functor $F_R$ is defined by
$$F^*M=R^{(1)}\otimes_RM$$
$$F^*(M\xrightarrow{\phi}N)=(R^{(1)}\otimes_RM\xrightarrow{id\otimes\phi}R^{(1)}\otimes_RN)$$
for all $R$-modules $M,N$ and all $R$-module homomorphisms $\phi$, where $F^*M$ acquires its $R$-module structure via the left $R$-module structure on $R^{(1)}$. Notice that $F^*$ is an exact functor by Kunz's theorem (Theorem 2.1 in \cite{kunz}).\par
If $M$ is a free $R$-module with an $R$-basis $\{b_1,b_2,\dots\}$, then the $R$-homomorphism $F^*M\to M$ given by
\begin{equation}
\label{isomorphism-free-module}
\sum_i r_i\otimes s_ib_i\mapsto \sum_i r_is^p_ib_i
\end{equation}
is an isomorphism.\par
For any $R$-module $M$, the isomorphism
\begin{equation}
\label{pullback-commute-Hom}
F^*(\Hom_R(M,R))\to \Hom_R(F^*M,R)
\end{equation}
is given by $r\otimes \phi\mapsto \psi$ where $\psi$ is defined by $\psi(s\otimes m)=rs\phi(m)^p$ for all $r,s\in R^{(1)},m\in M, \phi\in \Hom_R(M,R)$.\par
For any ideal $J$ in $R$, the homomorphism $F^*(R/J)\xrightarrow{r'\otimes \bar{r}\mapsto \overline{r'r^p}}R/J^{[p]}$ is an isomorphism.\par

Next we will discuss the grading on functors associated with $F$.

\begin{discussion}{\bf (grading on $F^*M$)}
\label{grading-on-pullback}
{\rm In order to make the isomorphism $F^*R\to R$ (\ref{isomorphism-free-module}) a degree-preserving $R$-homomorphism, the grading one can put on $F^*R$ is given by
$$\deg(r'\otimes r)=\deg(r')+p\deg(r)$$
where $r'\in R^{(1)}$ and $r\in R$ are homogeneous elements. In general, for any graded $R$-module $M$, the grading on $F^*M$ is given by
$$\deg(r'\otimes m)=\deg(r')+p\deg(m)$$
where $r'\in R^{(1)}$ and $m\in M$ are homogeneous elements. \par
It is easy to check that the isomorphism $F_R(R/J)\to R/J^{[p]}$ is degree-preserving for all ideals $J$.}
\end{discussion}

\begin{discussion}[grading on $F_*M$]
\label{grading-on-pushforward}
{\rm Let $M$ be a $\mathbb{Z}$-graded $R$-module. By abuse of notation, for any element $z\in M$, we denote by $F_*(z)$ the corresponding element in $F_*M$. Then we put a $\mathbb{Z}[\frac{1}{p}]$-grading on $F_*M$ as follows, 
$$\deg(F_*(z))=\frac{1}{p}\deg(z)$$
for all homogeneous elements $z\in M$. This grading respects the $R$-module structure on $F_*M$, i.e., for any homogeneous element $r\in R$, one has $\deg(rF_*(z))=\deg(F_*(r^pz))=\deg(r)+\deg(F_*(z))$.\par
Given any $\phi\in \Hom_R(F_*M,N)$ where $N$ is a $\mathbb{Z}$-graded $R$-module, if there is an $a\in \mathbb{Z}[\frac{1}{p}]$ such that $\deg(\phi(F_*(z)))=\deg(F_*(z))+a$ for all homogeneous $z$, we define the degree of $\phi$ to be $a$ and call $\phi$ homogeneous. We denote by ${^*\Hom}_R(F_*M,N)$ the submodule of $\Hom_R(F_*M,N)$ generated by homogeneous elements. It is clear that when $M$ is finitely generated ${^*\Hom}_R(F_*M,N)=\Hom_R(F_*M,N)$. When this is the case we will just write $\Hom_R(F_*M,N)$ assuming no confusion will arise. Hence if $M$ is a finitely generated $\mathbb{Z}$-graded $R$-module, $\Hom_R(F_*M,N)$ is a $\mathbb{Z}[\frac{1}{p}]$-graded $R$-module for all $\mathbb{Z}$-graded modules $N$. Consequently, $\Ext^t_R(F_*M,N)$ is also a $\mathbb{Z}[\frac{1}{p}]$-graded $R$-module for all integers $t$ and $\mathbb{Z}$-graded modules $N$ (since, if $\cdots\to K_1\to K_0\to M\to 0$ be a graded free resolution of $M$, then $\cdots\to F_*K_1\to F_*K_0\to F_*M\to 0$ is also a graded free resolution of $F_*M$).\par

This $\mathbb{Z}[\frac{1}{p}]$-grading has the same effect on degree-shifts, i.e., for any $\mathbb{Z}$-graded $R$-module $M$, one has $F_*(M(l))=(F_*M)(\frac{l}{p})$ for all integers $l$.
}
\end{discussion} 

\begin{discussion}
{\rm For any $\mathbb{Z}[\frac{1}{p}]$-graded $R$-module $M$, we denote by $\Int(M)$ the direct sum of graded pieces of $M$ with integer degrees (hence $\Int(M)$ is naturally a $\mathbb{Z}$-graded $R$-module). We observe that, for any $\mathbb{Z}$-graded $R$-modules $M$ and $N$ with $M$ finitely generated, we have 
$$\Int(\Hom_R(F_*M,N))=\Hom_R(\Int(F_*M),N)$$
and we reason as follows. Given any $\phi\in \Int(\Hom_R(F_*M,N))$ the image of any homogeneous element with fractional degree under $\phi$ must be 0, hence $\phi$ is naturally a homomorphism from $\Int(F_*M)$ to $N$. This proves that $\Int(\Hom_R(F_*M,N))\subseteq\Hom_R(\Int(F_*M),N)$. The other containment is clear. Consequently, we have
\[\Int(\Ext^t_R(F_*M,N))=\Ext^t_R(\Int(F_*M),N),\ {\rm for\ all\ integers\ }t,\]
and hence
\begin{equation}
\label{int-commute-EE}
\Int(\Ext^s_R(\Ext^t_R(F_*M,N),N'))=\Ext^s_R(\Ext^t_R(\Int(F_*M),N),N')
\end{equation}
for all integers $s,t$.}
\end{discussion}

\begin{discussion}
\label{equivalence}
{\rm
Set
\[\cA=\{{\rm graded\ }p{\rm -linear\ structures\ on\ }M\},\]
\[\cB=\{{\rm degree-preserving\ }R{\rm -linear\ maps\ }M\to F_*M\},\]
and
\[\cC=\{{\rm degree-preserving\ }R{\rm -linear\ maps\ }F^*M\to M\}.\]
We have natural bijections among $\cA$, $\cB$ and $\cC$ as follows.
\[\xymatrix{
\cA \ar@<1ex>[r]^{\sigma} & \cB \ar@<1ex>[r]^{\gamma} \ar@<1ex>[l]^{\sigma'} & \cC \ar@<1ex>[l]^{\gamma'}
}\]
Let $\epsilon:M\to F_*M$ be the identity map on the underlying additive group (we will also use $\epsilon$ to denote $F_*M\to M$ assuming no confusion will arise). Then $\sigma:\cA\to \cB$ is defined by $\sigma(f)=\epsilon\circ f$; and $\sigma':\cB\to \cA$ is defined by $\sigma'(g)=\epsilon\circ g$. Also $\gamma:\cB\to \cC$ is defined by $\gamma(g)(r\otimes m)=rg(m)$ for all $r\in R^{(1)}$ and $m\in M$; while $\gamma':\cC\to \cB$ is defined by $\gamma'(h)=h(1\otimes m)$ for all $m\in M$. It is straightforward to check that, $\sigma$ and $\sigma'$ ($\gamma$ and $\gamma'$ respectively) are inverses to each other.
}
\end{discussion}

Recall from Section 2 that each finite morphism $f:Y_1\to Y_2$ induces a functor $f^!$ as follows
\[f^!(-)=\bar{f}^*\underline{\underline{R}}\underline{\Hom}(f_*\mathcal{O}_{Y_1},-).\]
In particular, since the Frobenius $F:\Spec(R)\to \Spec(R)$ is a finite morphism, we have a functor $F^!$ (with $f$ replaced by $F$ in the above definition of $f^!$). This functor $F^!$ can actually be viewed as a functor from the category of $R$-modules to itself. 

\begin{discussion}[Grading on $F^!M$]
\label{grading-on-F-!}
{\rm To avoid any ambiguity arising from the different roles played by the two $R$s in $\Spec(R)\xrightarrow{F}\Spec(R)$, we will write the Frobenius endomorphism on $\Spec(R)$ as $F:\Spec(R')\to \Spec(R)$. Since $\Spec(R)$ is affine and $F$ is a finite flat morphism, the functor $F^!$ is actually a functor from the category of $R$-modules to the category of $R'$-modules (not just a functor between derived categories):
\[F^!N=\Hom_R(R',N)\ {\rm as\ an}\ R'{\rm -module}\]
for all $R$-modules $N$ (cf. \cite[page 171]{rd}).\par
The module $\Hom_R(R',N)$ can be viewed as an $R$-module and as an $R'$-module. Again to avoid any ambiguity arising from this, we will use ${\sideset{_R}{_R}\Hom} (R',N)$ to denote $\Hom_R(R',N)$ viewed as an $R$-module; and we will use ${\sideset{_{R'}}{_R}\Hom}(R',N)$ to denote $\Hom_R(R',N)$ viewed as an $R'$-module. This way, $F^!$ is actually given by
\[F^!(N)={\sideset{_{R'}}{_R}\Hom}(R',N)\]
for all $R$-modules $N$; and
\begin{equation}
\label{push-forward-F-!} 
F_*({\sideset{_{R'}}{_R}\Hom}(R',N))={\sideset{_R}{_R}\Hom}(R',N).
\end{equation}
If $N$ is graded, then the grading on ${\sideset{_R}{_R}\Hom}(R',N)$ is induced from the grading on $N$ and the grading on $F_*R$ according to Discussion \ref{grading-on-pushforward} ($R'$ as an $R$-module is nothing but $F_*R$); while the grading on ${\sideset{_{R'}}{_R}\Hom} (R',N)$ is obtained by multiplying the just defined degrees by $p$; i.e. an element $\phi\in \Hom_R(R',M)$ is homogeneous in ${\sideset{_R}{_R}\Hom}(R',M)$ if and only if it is homogeneous in ${\sideset{ _{R'}}{_R}\Hom}(R',M)$ and $\deg_{R'}(\phi)=p\cdot \deg_R(\phi)$. It is not hard to check that with these degrees ${\sideset{_R}{_R}\Hom}(R',M)$ is a graded $R$-module, i.e. $\deg_R(r\cdot\phi)=\deg(r)+\deg_R(\phi)$ and  ${\sideset{_{R'}}{_R}\Hom}(R',M)$ is a graded $R'$-module, i.e. $\deg_{R'}(r'\cdot\phi)=\deg(r')+\deg_{R'}(\phi)$ for all homogeneous elements $r\in R$, $r'\in R'$ and $\phi\in \Hom_R(R',M)$ and (\ref{grading-on-pushforward}) is degree-preserving. Finally, we point out that 
\begin{equation}
\label{degree-shift-F-!}
F^!(N(\ell))=(F^!N)(p\cdot \ell)
\end{equation}
 as graded $R'$-modules, because} 
\[F_*(F^!(N(\ell)))={\sideset{_R}{_R}\Hom}(R',N(\ell))={\sideset{_R}{_R}\Hom}(R',N)(\ell)=F_*((F^!N)(p\cdot \ell)).\]
\end{discussion}

Given our grading on $F^!M$, we have the following result on the duality for the Frobenius morphism.

\begin{pro}
\label{Frobenius-duality-degree-preserving}
For any finitely generated graded $R$-module $M$, the $R$-module isomorphism 
\[F_*\Ext^t_R(M,F^!\Omega)\xrightarrow{\sim} \Ext^t_R(F_*M,\Omega)\]
induced by Corollary \ref{duality-sheaf-ext} is degree-preserving.
\end{pro}
\begin{proof}[Proof]
As before we will write $F$ as $F:\Spec(R')\to \Spec(R)$. (Thus the Frobenius map now reads $R\to R'$.) Under this notation, we will have to prove that the $R$-module isomorphism
\[F_*\Ext^t_{R'}(M',F^!\Omega)\xrightarrow{\sim}\Ext^t_R(F_*M',\Omega)\]
is degree-preserving for all finitely generated graded $R'$-module $M'$.\par

It is straightforward and left to the reader to check that the following $R$-module isomorphism (which induces Grothendieck Duality for finite morphisms \cite[page 171]{rd}) is degree-preserving
\begin{align}
\Delta:{\sideset{_R}{_{R'}}\Hom}(M',{\sideset{_{R'}}{_R}\Hom}(R',N))&\xrightarrow{\sim}\Hom_R(F_*M',N),\notag\\
 f&\mapsto g,\ g(z')=f(z')(1)\notag
\end{align}
where $M'$ is a finite graded $R'$-module and $N$ is a graded $R$-module, $z'$ is an element in $M'$ and at the same time an element of $F_*M'$ (because $M' = F_*M'$ as sets). \par
Let $0\to H'_{n+1}\to \cdots\to\cdots\to H'_0\to M'\to 0$ be a graded free resolution of $M'$ over $R'$. Applying $\Delta$ to the case $M'=H'_{t}$ and $N=\Omega$, we see that the isomorphism 
\[F_*\Ext^t_{R'}(M',F^!\Omega)\xrightarrow{\sim}\Ext^t_R(F_*M',\Omega)\]
is degree-preserving as desired.
\end{proof}

We end this section with the following easy (and very likely well-known) lemma.
\begin{lem}
\label{canonical-module-preserved}
There is a degree-preserving isomorphism
\[F^!\Omega\cong\Omega.\]
\end{lem} 
\begin{proof}[Proof]
It is straightforward to check that the classical Cartier Isomorphism $F_*\Omega\to \Omega$ is degree-preserving and hence its adjoint\footnote{The adjoint $\Omega\xrightarrow{\sim} F^!\Omega$ is actually the isomorphism induced by (\ref{R-prime-linear-isomorphism}) with suitable degree-shifts. One can also prove directly that it is degree-preserving.} $\Omega\xrightarrow{\sim} F^!\Omega$ is also degree-preserving.
\end{proof}

\section{A $p$-linear structure on $\mathcal{M}^{i,j}_0$}
This section and next section are devoted to proving Theorem \ref{indp-f-M}, whose proof is long and technical. Before we proceed to the proof, we explain the idea behind the proof. We need to construct a graded $p$-linear structure $\varphi$ on $\cM^{i,j}$ with two main properties: on one hand, $\varphi_0$ should be independent of the embedding; on the other hand, it should be compatible with Theorem \ref{stable-M-i-j} (see Remark \ref{compatibility} and Proposition \ref{limit-lambda} for details). We will achieve this in, roughly speaking, 3 steps. First, we construct two different looking graded $p$-linear structures $\varphi$ and $\varphi'$ on $\cE^{i,j}(R/I)$: one is induced by $R/I\xrightarrow{\hat{r}\mapsto \hat{r}^p} R/I$ while the other comes from Corollary \ref{f-finite-gamma}. Then we will prove that $\varphi$ and $\varphi'$ coincide. All of this is accomplished in this section. Finally, in next section we prove that $\varphi_0$ actually is independent of the embedding; this will be the most technical step.\par 

We begin with two lemmas which assert that both $F^*$ and $F_*$ commute with $\cE^{i,j}(-)$ in the graded sense.

\begin{lem}
\label{pullback-commute-cE}
There is a degree-preserving isomorphism:
\[F^*\cE^{i,j}(M)\xrightarrow{\sim}\cE^{i,j}(F^*M).\]
\end{lem}
\begin{proof}[Proof]
With our grading on $F^*$, it is easy to check that, the isomorphism (\ref{pullback-commute-Hom}) is degree-preserving, hence, so are induced isomorphisms
$$F^*(\Ext^s_R(M,R))\to \Ext^s_R(F^*M,R)$$
for all integers $s$. 
Therefore the induced isomorphisms
\begin{align}
F^*\cE^{i,j}(M)&=F^*(\Ext^{n+1-i}_R(\Ext^{n+1-j}(M,R),R))\notag\\ &\to \Ext^{n+1-i}_R(F^*(\Ext^{n+1-j}_R(M,R)),R)\notag \\
                            &\to \Ext^{n+1-i}_R(\Ext^{n+1-j}_R(F^*M,R),R)\notag\\
                            &=\cE^{i,j}(F^*M)\notag
\end{align}
are degree-preserving for all $i,j$.
\end{proof}

Over $\Spec(R)$, it is clear that 
\[\E^{-i,-j}(\Spec(R);M)=\cE^{i,j}(M)\]
for all $R$-modules $M$, since $\omega^{\bullet}_{\Spec(R)}=R[-n-1]$ the complex with $R$ in degree $n+1$ and 0 elsewhere. To simplify the notations, we will write $\E^{-i,-j}(M)$, with $\Spec(R)$ omitted, for all $R$-modules $M$. One may view $\E^{-i,-j}(-)$ as a covariant functor from the category of $R$-modules to itself. If $M$ is a graded $R$-module, the gradings on $M$ and $R$ induce a natural grading on $\E^{s,t}(M)$.
\begin{lem}
\label{pushforward-commute-cE}
The following isomorphism induced from Corollary \ref{f-finite-gamma} is degree-preserving
$$\cE^{i,j}(F_*M)=\E^{-i,-j}(F_*M)\xrightarrow{\sim} F_*\E^{-i,-j}(M)=F_*\cE^{i,j}(M).$$
\end{lem}
\begin{proof}[Proof]
It follows directly from Proposition \ref{Frobenius-duality-degree-preserving} and Lemma \ref{canonical-module-preserved}.
\end{proof}

\begin{discussion}[Graded $p$-linear structures $\varphi$ and $\varphi'$ on $\cE^{i,j}(R/I)$]
\label{graded-p-structures-on-cE}
{\rm 
First, We construct $\varphi$:
\[\varphi:\cE^{i,j}(R/I)\xrightarrow{\sigma_1}F^*\cE^{i,j}(R/I)\xrightarrow{\sigma_2}\cE^{i,j}(F^*(R/I))\xrightarrow{\sigma_3}\cE^{i,j}(R/I)\]
where $\sigma_1$ is defined by $z\mapsto 1\otimes z$; while $\sigma_2$ comes from Lemma \ref{pullback-commute-cE}; and $\sigma_3$ is induced from $F^*(R/I)\cong R/I^{[p]}\twoheadrightarrow R/I$. It is clear that $\varphi$ is a graded $p$-linear structure on $\cE^{i,j}(R/I)$.\par
Next, we construct $\varphi'$:
\[\varphi':\cE^{i,j}(R/I)\xrightarrow{\tau_1}\cE^{i,j}(F_*(R/I))\xrightarrow{\tau_2}F_*\cE^{i,j}(R/I)\xrightarrow{\tau_3}\cE^{i,j}(R/I),\]
where $\tau_1$ is induced by the degree-preserving $R$-linear map $R/I\to F_*(R/I)$; and $\tau_2$ is the isomorphism in Lemma \ref{pushforward-commute-cE}; and $\tau_3$ is the identity map on the underlying additive group. It is clear that $\varphi'$ is also a graded $p$-linear structure on $\cE^{i,j}(R/I)$.\par
}
\end{discussion}

As we will see that $\varphi$ and $\varphi'$ coincide. The reason for us to introduce two such graded $p$-linear structures is that, $\varphi$ makes the transition maps in the inverse systems 
$$\{\cE^{i,j}(R/I^{[p^e]})\}_e \ {\rm and}\ \{F^{e*}(\cE^{i,j}(R/I))\}_e$$
compatible (see Remark \ref{compatibility} for details); while it is easier to work with $\varphi'$ when proving the embedding-independence of $\varphi_0$.

From the above discussion, it is clear that 
\begin{equation}
\label{F-commute-E-s-t}
F_*\E^{s,t}(-)\cong \E^{s,t}(F_*(-))\ {\rm and}\ F^*\E^{s,t}(-)\cong \E^{s,t}(F^*(-)) 
\end{equation}

The rest of this section is devoted to proving that $\varphi$ and $\varphi'$ coincide. To avoid any ambiguity arising from the different roles played by the two $R$s in $\Spec(R)\xrightarrow{G}\Spec(R)$, we will write $F:\Spec(R')\to \Spec(R)$, i.e., we will denote by $R'$ the target ring in $R\xrightarrow{r\mapsto r^p}R$. If $M'$ is an $R'$-module, then $F_*M'$ is the same abelian group as $M'$ but the $R$-module structure is by $r\cdot m=r^pm$ for all $m'\in M'$. It is well-known that $F^*$ is left-adjoint to $F_*$, i.e., for all $R$-modules $M$ and all $R'$-modules $N'$, one has the bijiective adjunction map
\begin{equation}
\label{F^*-left-adjoint-F_*}
\Theta_{M,N'}:\Hom_{R'}(F^*M,N')\xrightarrow{\sim} \Hom_R(M,F_*N') 
\end{equation}
where $\Theta_{M,N'}$ is given by $\phi\mapsto (M\xrightarrow{m\mapsto 1\otimes m} F_*F^*M\xrightarrow{F_*\phi}F_*N')$. In particular, one also has a bijection 
\[ \Theta_{\E^{s,t}(M),\E^{s,t}(N')}:\Hom_R(F^*\E^{s,t}(M),\E^{s,t}(N'))\to\Hom_R(\E^{s,t}(M),F_*\E^{s,t}(N')).\]
Furthermore, it is shown in \cite{lyzhzh} that $F^*$ is right-adjoint to $F_*$ as well, i.e., one has a bijection (actually an $R'$-linear isomorphism)
\begin{equation}
\label{F^*-right-adjoint-F_*}
\Lambda_{N',M}:\Hom_{R'}(N',F^*M)\xrightarrow{\sim}\Hom_R(F_*N',M)
\end{equation}
for all $R'$-modules $N'$ and $R$-modules $M$, where $\Lambda_{N',M}$ is defined as follows. It is clear that $R'$ is a free $R$-module of rank $p^{n+1}$ with a basis $\{b_{\bar{e}}=x^{e_0}_0\cdots x^{e_n}_n|0\leq e_i\leq p-1\}$ and that $\Hom_R(R',R)$ is also a free $R$-module of rank $p^{n+1}$ with a dual basis $\{\tilde{b}_{\bar{e}}|0\leq e_i\leq p-1\}$, i.e., $\tilde{b}_{\bar{e}}(b_{\bar{e'}})=1$ if $\bar{e}=\bar{e'}$ and $\tilde{b}_{\bar{e}}(b_{\bar{e'}})=0$ if $\bar{e}\neq\bar{e'}$, where $\bar{e}$ (respectively, $\bar{e'}$) denotes $(e_0,\dots,e_n)$ (respectively, $(e'_0,\dots,e'_n)$). For all $R'$-modules $N'$ and $R$-modules $M$, each $R'$-linear map $\psi:N'\to F^*M$ has the form $\psi=\oplus_{\bar{e}}(b_{\bar{e}}\otimes_R\psi_{\bar{e}})$ where $\psi_{\bar{e}}:N'\to M$ are $R$-linear. Then  
\[\Lambda_{N',M}(\psi=\oplus_{\bar{e}}(b_{\bar{e}}\otimes_R\psi_{\bar{e}}))=\psi_{\overline{p-1}},\]
where $\overline{p-1}$ denotes $(p-1,\dots,p-1)$. It is shown in \cite[\S4]{lyzhzh} that $\Lambda_{N',M}$ is an $R'$-linear isomorphism. Note that, given the basis $\{b_{\bar{e}}\}$ of $R'$ as a free $R$-module and its dual basis $\{\tilde{b}_{\bar{e}}\}$, one has an $R'$-linear isomorphism 
\begin{equation}
\label{R-prime-linear-isomorphism}
R'\xrightarrow{\sim} \Hom_R(R',R),\ b_{\bar{e}}\mapsto \tilde{b}_{\overline{p-1}-\bar{e}}\ .
\end{equation}

It is clear from the discussion at the beginning of this section that we have an isomorphism
\begin{equation}
\label{F-upper-star-commute-with-Ext}
\Ext^{t}_{R'}(F^*M,R')\cong F^*\Ext^{t}_R(M,R)
\end{equation}
for all $R$-modules $M$ and integers $t$. Combining Corollary \ref{duality-sheaf-ext} and (\ref{R-prime-linear-isomorphism}), we have an isomorphism
\begin{equation}
\label{F-lower-star-commute-with-Ext}
\Ext^{t}_R(F_*N',R)\cong F_*\Ext^{t}_{R'}(N',R')
\end{equation}
for all $R'$-modules $N'$ and integers $t$. Then one has the following commutative diagrams (see \cite[\S4]{lyzhzh} for details):
{\tiny
\begin{equation}
\label{theta-lambda}
\xymatrix{
\Hom_{R'}(F^*M,N')\ar[r]^{\Theta_{M,N'}} \ar[d]  &  \Hom_R(M,F_*N')\ar[d] \\
\Hom_{R'}(\Ext^{t}_{R'}(N',R'), \Ext^{t}_{R'}(F^*M,R'))\ar[d]^{y_1}_{\sim} &  \Hom_R(\Ext^{t}_R(F_*N',R), \Ext^{t}_R(M,R))\ar[d]^{z_1}_{\sim}\\
\Hom_{R'}(\Ext^{t}_{R'}(N',R'), F^*\Ext^{t}_R(M,R)) \ar[r]^{\sim} & \Hom_R(F_*\Ext^{t}_{R'}(N',R'), \Ext^{t}_R(M,R))
}
\end{equation}}
and 
{\tiny
\begin{equation}
\label{lambda-theta}
\xymatrix{
\Hom_{R'}(N',F^*M)\ar[r]^{\Lambda_{N',M}} \ar[d]  &  \Hom_R(F_*N',M)\ar[d] \\
\Hom_{R'}(\Ext^{s}_{R'}(F^*M,R'), \Ext^{s}_{R'}(N',R'))\ar[d]^{y_2}_{\sim} &  \Hom_R(\Ext^{s}_R(M,R), \Ext^{s}_R(F_*N',R))\ar[d]^{z_2}_{\sim}\\
\Hom_{R'}(F^*\Ext^{s}_R(M,R), \Ext^{s}_{R'}(N',R')) \ar[r]^{\sim} &  \Hom_R(\Ext^{s}_R(M,R), F_*\Ext^{s}_{R'}(N',R'))
}
\end{equation}}
where the top vertical maps in both diagrams are induced by functoriality while $y_1$ and $y_2$ (resp. $z_1$ 
and $z_2$) are induced by (\ref{F-upper-star-commute-with-Ext}) (resp. (\ref{F-lower-star-commute-with-Ext})).\par

Combining (\ref{theta-lambda}) and (\ref{lambda-theta}), one has immediately the following lemma which draws connections between $\Theta_{M,N'}$ and $\Theta_{\E^{s,t}(M),\E^{s,t}(N')}$.
\begin{lem}
\label{adjunction}
Given $\alpha\in\Hom_{R'}(F^*M, N')$ and $\beta\in \Hom_R(M,F_*N')$, let $\widetilde{\alpha}$ and $\widetilde{\beta}$ be the following induced maps $$\widetilde{\alpha}:F^*\E^{s,t}(M)\stackrel{(\ref{F-commute-E-s-t})}{\cong} \E^{s,t}(F^*M)\xrightarrow{\E^{s,t}(\alpha)}\E^{s,t}(N')$$
and$$\widetilde{\beta}:\E^{s,t}(M)\xrightarrow{\E^{s,t}(\beta)}\E^{s,t}(F_*N')\stackrel{(\ref{F-commute-E-s-t})}{\cong} F_*\E^{s,t}(N').$$
If $\beta=\Theta_{M,N'}(\alpha)$, then $\widetilde{\beta}=\Theta_{\E^{s,t}(M),\E^{s,t}(N')}(\widetilde{\alpha})$. In other words, we have a commutative diagram
\[\xymatrix{
\Hom_{R'}(F^*M,N') \ar[rr]^{\Theta_{M,N'}}_{\sim} \ar[d]^{\epsilon_1} & & \Hom_R(M,F_*N') \ar[d]^{\epsilon_2} \\
\Hom_{R'}(\E^{s,t}(F^*M),\E^{s,t}(N')) \ar[d]^{d_1}_{\sim} & & \Hom_R(\E^{s,t}(M),\E^{s,t}(F_*N')) \ar[d]^{d_2}_{\sim}\\
\Hom_{R'}(F^*\E^{s,t}(M),\E^{s,t}(N')) \ar[rr]^{\Theta_{\E^{s,t}(M),\E^{s,t}(N')}}_{\sim} & & \Hom_R(\E^{s,t}(M),F_*\E^{s,t}(N'))
}\]
where $\epsilon_1$ and $\epsilon_2$ are induced from the fact that $\E^{s,t}(-)$ is a covariant functor; $d_1$ and $d_2$ are induced by (\ref{F-commute-E-s-t}).
\end{lem}

We are ready to show that $\varphi$ and $\varphi'$ coincide.

\begin{pro}
\label{frobenius-M-i-j}
The graded $p$-linear structures $\varphi$ and $\varphi'$ on $\cE^{i,j}(R/I)$ coincide.
\end{pro}
\begin{proof}[Proof]
Let the map $\alpha:F^*(R/I)\to R/I$ be induced by the natural surjection $R/I^{[p]}\to R/I$, then it is clear that $\beta=\Theta_{R/I,R/I}(\alpha):R/I\to F_*(R/I)$ is the map $R/I\xrightarrow{r\mapsto r^p} F_*(R/I)$. By Lemma \ref{adjunction}, we know that 
\begin{equation}
\label{adjunction-equality}
\Theta_{\E^{-i,-j}(R/I),\E^{-i,-j}(R/I)}\circ d_1\circ \epsilon_1(\alpha)= d_2\circ\epsilon_2\circ\Theta_{R/I,R/I}(\alpha)
\end{equation}
In terms of $F_*$ and $F^*$, it is clear that $\varphi$ is the composition of
\begin{align} 
\cE^{i,j}(R/I)=\E^{-i,-j}(R/I) &\to F_*F^*\E^{-i,-j}(R/I)\notag\\ & \xrightarrow{\sim}F_*\E^{-i,-j}(F^*(R/I))\notag\\
 & \to F_*\E^{-i,-j}(R/I)=F_*\cM^{i,j}\notag
\end{align}
which is the left hand side of (\ref{adjunction-equality}) with $M=N'=R/I$. On the other hand, the graded $p$-linear structure $\varphi'$ is induced by $\beta:R/I\to F_*(R/I)$, i.e., $\varphi'$ is the right hand side of (\ref{adjunction-equality}) with $M=N'=R/I$. Therefore $\varphi$ and $\varphi'$ coincide.
\end{proof}


\section{$\varphi_0$ is independent of the embedding}
In this section we will finish the last step of the proof of Theorem \ref{indp-f-M}, i.e., $\varphi_0$ is independent of the embedding.\par

For all schemes over $k$, throughout this section we always assume that the Frobenius endomorphism is a finite morphism (for example, any scheme of finite type over $k$ has a finite Frobenius endomorphism). To ease the notations, we will write $\E^{*,*}(M)$ instead of $\E^{*,*}(\Spec(R),M)$ for all $R$-modules $M$.\par

Recall that in the construction of $\varphi'$ (Discussion \ref{graded-p-structures-on-cE}) we define an $R$-linear isomorphism
\[
\tau_2:\cE^{i,j}(F_*(R/I))\cong F_*\cE^{i,j}(R/I)
\]
which is degree-preserving by Lemma \ref{pushforward-commute-cE}. Hence we have a $k$-linear isomorphism 
\[(\tau_2)_0:\cE^{i,j}(F_*(R/I))_0\xrightarrow{\sim} F_*\cE^{i,j}(R/I)_0.\]

By virtue of (\ref{int-commute-EE}), we can view $(\tau_2)_0$ as
\[(\tau_2)_0:\E^{-i,-j}(\Int(F_*(R/I)))_0\xrightarrow{\sim}(F_*\E^{-i,-j}(R/I))_0.\]

\begin{pro}
\label{independent-tau-2}
The isomorphism
\[(\tau_2)_0:\E^{-i,-j}(\Int(F_*(R/I)))_0\xrightarrow{\sim}(F_*\E^{-i,-j}(R/I))_0\]
depends only on $\cO_X$.
\end{pro}

Before we proceed to the proof of this proposition, we need the following remark and lemma.

\begin{rem}
\label{commute-G-F-k}
{\rm Let $R$ and $F$ be as above and let $F_k$ denote the Frobenius endomorphism on $\Spec(k)$. Then
\begin{enumerate}
\item For any $\mathbb{Z}$-graded $R$-module $M$, one has $(F_*M)_0=F_{k*}(M_0)$ since $(F_*M)_0$ is the same as $M_0$ as abelian groups and the $k$ vector space structure on $(F_*M)_0$ is via $k\xrightarrow{c\mapsto c^p}k$. And hence $(\Int(F_*M))_0=F_{k*}(M_0)$.
\item Let $F_{\mathbb{P}^n}$ denote the Frobenius endomorphism on $\mathbb{P}^n$. One has that $F^{\Gamma}_{\mathbb{P}^n*}=F_{k*}$, since $\Gamma(\mathbb{P}^n,\cO_{\mathbb{P}^n})=k$. 
\item Given the degree on $F_*M$ for a $\mathbb{Z}$-graded $R$-module $M$ as in Discussion \ref{grading-on-pushforward}, then it is well-known that $\widetilde{\Int(F_*M)}=F_{\mathbb{P}^n_*}\widetilde{M}$, where $\widetilde{\Int(F_*M)}$ and $\widetilde{M}$ are sheaves associated to $\Int(F_*M)$ and $M$ respectively on $\mathbb{P}^n$. For completeness we give a proof as follows. It suffices to prove that, for integers $l\gg 0$, we have 
\[H^0(\mathbb{P}^n,(F_{\mathbb{P}^n*}\widetilde{M})(l))=(\Int(F_*M))_l,\]
but it is clear that
\begin{align}
&\quad H^0(\mathbb{P}^n,(F_{\mathbb{P}^n*}\widetilde{M})(l))\notag\\
&=H^0(\mathbb{P}^n,F_{\mathbb{P}^n*}\widetilde{M}\otimes\cO_{\mathbb{P}^n}(l))\notag\\
&=H^0(\mathbb{P}^n,F_{\mathbb{P}^n*}(\widetilde{M}\otimes F^*_{\mathbb{P}^n}(\cO_{\mathbb{P}^n}(l)))\tag{\dag}\\
&=H^0(\mathbb{P}^n,F_{\mathbb{P}^n*}(\widetilde{M}\otimes\cO_{\mathbb{P}^n}(pl))\notag\\
&=H^0(\mathbb{P}^n,F_{\mathbb{P}^n*}(\widetilde{M}(pl)))\notag\\
&=H^0(\mathbb{P}^n,\widetilde{M}(pl))\ {\rm since}\ F\ {\rm is\ an\ affine\ morphism}\notag\\
&\stackrel{l\gg 0}{=}M_{pl}\notag\\
&=\Int(F_*M))_l,\notag
 \end{align}
where ($\dag$) follows from the Projection Formula (\cite[page 124]{ag}).
\end{enumerate}}
\end{rem}

\begin{lem}
\label{commute-D-F-k}
Let $R$, $\Omega$ be as above, let $\D(M)$ denote the graded Matlis dual of any graded $R$-module $M$ as in {\rm (GLD)} and let $\D(V)$ denote the dual space of any $k$-vector space $V$ (i.e., $\D(V)=\Hom_k(V,k)$). Then
\begin{enumerate}
\item[(a)] For any finite dimensional $k$ vector space $V$, one has a functorial isomorphism
\[\D(F_{k*}V)\cong F_{k*}D(V).\]
\item[(b)] There is a functorial $k$-linear isomorphism 
\[\D((F_*M)_0)\xrightarrow{\sim}(F_*(\D(M)))_0=F_{k*}(\D(M)_0)=F_{k*}(\D(M_0))\]
for any $\mathbb{Z}$-graded $R$-module $M$.
\item[(c)] There is a functorial $k$-linear isomorphism
\[(F_*H^t_{\fm}(M))_0=F_{k*}(H^t_{\fm}(M)_0)\xrightarrow{\sim}H^t_{\fm}(\Int(F_*M))_0 .\]
for any $\mathbb{Z}$-graded $R$-module $M$.
\item[(d)] The map (\ref{degree-0-ext-sheaf-ext}) is compatible with Frobenius maps in the sense that, for any finitely generated $\mathbb{Z}$-graded $R$-module $M$, there is a commutative diagram (for $i\geq 1$) 

\begin{equation}
\label{diagram-degree-0-ext-sheaf}
\xymatrix{
\Ext^{n+1-i}_R(\Int(F_*M),\Omega)_0 \ar[r]^{u} \ar[d]^{v}_{\sim} & \Ext^{1-i}_{\mathbb{P}^n}(F_{\mathbb{P}^n*}\widetilde{M},\dualp)\ar[d]^{v'}_{\sim}\\
(\Int(F_*\Ext^{n+1-i}_R(M,\Omega)))_0 \ar@{=}[d] & F^{\Gamma}_{\mathbb{P}^n*}\Ext^{1-i}(\widetilde{M},\dualp)\ar@{=}[d] \\
F_{k*}(\Ext^{n+1-i}_R(M,\Omega)_0) \ar[r]^{u'=F_{k*}(u)} & F_{k*} \Ext^{1-i}_{\mathbb{P}^n}(\widetilde{M},\dualp)
}
\end{equation}
where $u$ and $u'$ are induced by (\ref{degree-0-ext-sheaf-ext}); $v$ is from Proposition \ref{Frobenius-duality-degree-preserving} and $v'$ is from Corollary \ref{duality-sheaf-ext}; finally the two vertical equalities follow from Discussion \ref{grading-on-pushforward}, Remark \ref{commute-G-F-k}(1) and the fact that $F^{\Gamma}_{\mathbb{P}^n*}=F_{k*}$.
\end{enumerate}
\end{lem}
\begin{proof}[Proof]
(a). This is clear.\par
(b). This holds because
\begin{align}
\D((F_*M)_0) &=\D(F_{k*}(M_0))=\Hom_k(F_{k*}(M_0),k)\notag\\
&\xrightarrow{\sim}F_{k*}\D(M_0)\notag\\
&=F_{k*}(\D(M))_0\notag\\
&=(F_*\D(M))_0\notag
\end{align}
and $\Hom_k(F_{k*}(M_0),k)\cong F_{k*}\D(M_0)$ is functorial by part (a).\par
(c). The first equality follows from Remark \ref{commute-G-F-k}(1). The isomorphism is given by 
\begin{align}
F_{k*}(H^t_{\fm}(M)_0)&\xrightarrow{\sim} F_{k*}(\D(\Ext^{n+1-t}_R(M,\Omega)_0))\ {\rm by\ (GLD)}\notag\\
&\xrightarrow{\sim}\D(F_{k*}(\Ext^{n+1-t}_R(M,\Omega)_0))\ {\rm by\ part\ (a)}\notag\\
&\xrightarrow{\sim}\D((\Int(F_*\Ext^{n+1-t}_R(M,\Omega)))_0)\ {\rm by\ Remark\ \ref{commute-G-F-k}}\notag\\
&\xrightarrow{\sim}\D(\Ext^{n+1-t}_R(\Int(F_*M),\Omega)_0)\ {\rm by\ Proposition\ \ref{canonical-module-preserved}}\notag\\
&\xrightarrow{\sim}H^t_{\fm}(\Int(F_*M))_0\ {\rm by\ (GLD)}\notag
\end{align}
which is functorial since all isomorphisms involved are functorial.\par
(d). In order to prove the commutativity of (\ref{diagram-degree-0-ext-sheaf}), we can refine it as follows

\begin{equation}
\label{refine-diagram}
\xymatrix{
\Ext^{n+1-i}_R(\Int(F_*M),\Omega)_0 \ar[r]^{v}_{\sim}\ar[d]^{u_1}_{\sim} & F_{k*}(\Ext^{n+1-i}_R(M,\Omega)_0)\ar[d]^{u'_1}\\
\D(H^i_{\fm}(\Int(F_*M))_0)\ar[r]^{v_1}_{\sim}\ar[d]^{u_2} & F_{k*}\D(H^i_{\fm}(M)_0) \ar[d]^{u'_2}\\
\D(H^{i-1}(\mathbb{P}^n,F_{\mathbb{P}^n*}\widetilde{M})) \ar[r]^{v_2}_{\sim}\ar[d]^{u_3}_{\sim} & F_{k*}\D(H^{i-1}(\mathbb{P}^n,\widetilde{M})) \ar[d]^{u'_3}\\
\Ext^{1-i}_{\mathbb{P}^n}(F_{\mathbb{P}^n*}\widetilde{M},\dualp) \ar[r]^{v'}_{\sim} & F_{k*}\Ext^{1-i}_{\mathbb{P}^n}(\widetilde{M},\dualp)
}
\end{equation}
where $u_3\circ u_2\circ u_1=u$ and $u'_3\circ u'_2\circ u'_1=u'$; more explicitly, $u_1$ and $u'_1$ are induced from (GLD), $u_2$ and $u'_2$ are induced from (LCSC) which are isomorphims when $i\geq 2$; $u_3$ and $u'_3$ are from Serre Duality; the map $v_1$ is given by 
\[\D(H^i_{\fm}(\Int(F_*M))_0)\xrightarrow[\sim]{{\rm by\ (c)}}\D(F_{k*}(H^i_{\fm}(M)_0))\xrightarrow[\sim]{{\rm by\ (a)}}F_{k*}(\D(H^i_{\fm}(M)_0));\]
and finally the map $v_2$ is given by
\[\D(H^{i-1}(\mathbb{P}^n,F_{\mathbb{P}^n*}\widetilde{M}))\cong \D(F_{k*}(H^{i-1}(\mathbb{P}^n,\widetilde{M})))\xrightarrow[\sim]{{\rm by\ (a)}}F_{k*}\D(H^{i-1}(\mathbb{P}^n,\widetilde{M})).\]
Now the commutativity of diagram (\ref{refine-diagram}) is clear because of the construction of the maps involved. \par
This completes the proof.
\end{proof}

We are now in the position to prove Proposition \ref{independent-tau-2}.
\begin{proof}[Proof of Proposition \ref{independent-tau-2}]First, we claim that there is a commutative diagram
{\tiny
\begin{equation}
\label{main-diagram-66}
\xymatrix{
\E^{-i,-j}(\Int(F_*(R/I)))_0 \ar[r]^{\zeta_3}\ar[d]^{q_1} & (F_*\E^{-i,-j}(R/I))_0=F_{k*}(\E^{-i,-j}(R/I)_0)\ar[d]^{q'_1}\\
\Ext^{1-i}_{\mathbb{P}^n}(\underline{\Ext}^{1-j}_{\mathbb{P}^n}(F_{\mathbb{P}^n*}\eta_*\cO_X,\dualp),\dualp)\ar[r]^{\zeta_2}_{\sim}\ar[d]^{q_2}_{\sim} & F_{k*}\Ext^{1-i}_{\mathbb{P}^n}(\underline{\Ext}^{1-j}_{\mathbb{P}^n}(\eta_*\cO_X,\dualp),\dualp)\ar[d]^{q'_2}_{\sim}\\
\E^{1-i,1-j}(X,F_{X*}\cO_X) \ar[r]^{\zeta_1}_{\sim} & F^{\Gamma}_{X*}\E^{1-i,1-j}(X,\cO_X)
}
\end{equation}
}
where $q_1$ and $q'_1$ are isomorphisms when $i\geq 2$; the maps $q_2$ and $q'_2$ are induced from Corollary \ref{f-finite-gamma} applied to $\eta:X\hookrightarrow\mathbb{P}^n$; the isomorphism $\zeta_2$ is given by
\begin{align}
\Ext^{1-i}_{\mathbb{P}^n}(\underline{\Ext}^{1-j}_{\mathbb{P}^n}(F_{\mathbb{P}^n*}\eta_*\cO_X,\dualp),\dualp)&\cong \Ext^{1-i}_{\mathbb{P}^n}(F_{\mathbb{P}^n*}\underline{\Ext}^{1-j}_{\mathbb{P}^n}(\eta_*\cO_X,\dualp),\dualp)\notag\\
&\cong F_{k*}\Ext^{1-i}_{\mathbb{P}^n}(\underline{\Ext}^{1-j}_{\mathbb{P}^n}(\eta_*\cO_X,\dualp),\dualp)\notag
\end{align}
where the first isomorphism is induced by the isomorphism \[\underline{\Ext}^{1-j}_{\mathbb{P}^n}(F_{\mathbb{P}^n*}\eta_*\cO_X,\dualp)\cong F_{\mathbb{P}^n*}\underline{\Ext}^{1-j}_{\mathbb{P}^n}(\eta_*\cO_X,\dualp)\] and the second isomorphism is the isomorphism $v'$ in diagram (\ref{refine-diagram}) with $\widetilde{M}=\underline{\Ext}^{1-j}_{\mathbb{P}^n}(\eta_*\cO_X,\dualp)$;  
and $\zeta_1$ is induced by Corollary \ref{f-finite-gamma} applied to $X\xrightarrow{F_X}X$ and is independent of the embedding. And we reason as follows. 
Let $\Delta_1$ be the diagram obtained from the diagram in Lemma \ref{lemma-degree-0-ext-sheaf} with $M=\Ext^{n+1-j}_R(\Int(F_*(R/I)),\Omega)$, $N=\Int(F_*(\Ext^{n+1-j}_R(R/I,\Omega)))$ and $M\to N$ being the degree-preserving isomorphism $M\cong N$ due to Proposition \ref{Frobenius-duality-degree-preserving}. (Note that 
\[\widetilde{M}=\underline{\Ext}^{1-j}(F_{\mathbb{P}^n*}\eta_*\cO_X,\dualp)\ {\rm and}\ \widetilde{N}=F_{\mathbb{P}^n*}\underline{\Ext}^{1-j}(\eta_*\cO_X,\dualp)\]
due to Proposition \ref{sheaf-ext} and Remark \ref{commute-G-F-k}(3).) Let $\Delta_2$ be the diagram obtained from Lemma \ref{commute-D-F-k}(d) by setting $M=\Ext^{n+1-j}_R(R/I,\Omega)$. (Note also that in this case $\widetilde{M}=\underline{\Ext}^{1-j}(\eta_*\cO_X,\dualp)$ due to Proposition \ref{sheaf-ext}.) The right column of $\Delta_1$ coincides with the left column of $\Delta_2$, so composing the horizontal maps produces the top square of (\ref{main-diagram-66}). Consequently, the top square of (\ref{main-diagram-66}) is commutative and $q_1$ and $q'_1$ are isomorphisms when $i\geq 2$. It is clear that the bottom square in diagram (\ref{main-diagram-66}) is commutative (this follows directly from the functoriality of Corollary \ref{f-finite-gamma}). This proves our claim.\par
To prove our proposition, we will consider two cases: $i\geq 2$ and $i=0,1$.\par
When $i\geq 2$, as we have seen that the maps $q_1$ and $q'_1$ in the diagram (\ref{main-diagram-66}) are isomorphisms, therefore $\zeta_3$ and $\zeta_2$ are independent of the embedding since $\zeta_1$ is so.\par
Next we consider the case $0\leq i\leq 1$. \par
We have an isomorphism
$$\rho:H^j_{\fm}(\Int(F_*(R/I)))_0\xrightarrow{\sim} F_{k*}H^j_{\fm}(R/I)_0$$
from Lemma \ref{commute-D-F-k}(c). From Proposition \ref{M-ij}(2), we know that $H^j_{\fm}(\Int(F_*(R/I)))_0$ (respectively, $H^j_{\fm}(R/I)_0$) only depends on $F_{X*}\cO_X$ (respectively, $\cO_X$); and it is straightforward (by dualizing the bottom two squares of the commutative diagram (\ref{refine-diagram})) to check that $\rho$ depends only on $X$.\par 
To finish the proof of our proposition, it suffices to prove the commutativity of the following diagram (since the columns are exact and both $\zeta_2$ and $\rho$ depend only on $X$)
 
\begin{equation}
\label{reduced-diagram}
\xymatrix{
0\ar[d] & 0\ar[d]\\
\E^{-1,-j}(\Int(F_*(R/I)))_0\ar[d]\ar[r]^{\zeta_3}_{\sim} &F_{k*}\E^{-1,-j}(R/I))_0\ar[d]\\
\E^{0,1-j}(\mathbb{P}^n;F_{\mathbb{P}^n*}\eta_*\cO_X)\ar[r]^{\zeta_2}_{\sim}\ar[d] &F_{k*}\E^{0,1-j}(\mathbb{P}^n;\eta_*\cO_X)\ar[d]\\
H^j_{\fm}(\Int(F_*(R/I)))_0\ar[d]\ar[r]^{\rho}_{\sim} & F_{k*}H^j_{\fm}(R/I)_0\ar[d]\\
\E^{0,-j}(\Int(F_*(R/I)))_0\ar[r]^{\zeta_3}_{\sim}\ar[d] & F_{k*}\E^{0,-j}(R/I)_0\ar[d]\\
0 & 0
}
\end{equation}

where the left column is obtained from (\ref{L-ij-0-1}) with $L=\Int(F_*(R/I))$ and the right column is obtained by applying $F^{\Gamma}_{\mathbb{P}^n*}=F_{k*}$ to (\ref{L-ij-0-1}) with $L=R/I$, and $\zeta_2$ is from Corollary \ref{duality-sheaf-ext} and is independent of the embedding.\par
The top square is exactly the diagram (\ref{main-diagram-66}) when $i=1$ and hence is commutative.\par
Next, we prove that the square in the middle is commutative. Set \newline $T=\Ext^{n+1-j}_R(R/I,\Omega)$ and $S=\Ext^{n+1-j}_R(\Int(F_*(R/I)),\Omega)$. Then we have a degree-preserving isomorphism
\[\Int(F_*(T))\cong \Int(\Ext^{n+1-j}_R(F_*(R/I),\Omega))=\Ext^{n+1-j}_R(\Int(F_*(R/I)),\Omega)=S\]
where the first isomorphism follows from Proposition \ref{Frobenius-duality-degree-preserving} and the second equality follows from Discussion \ref{grading-on-pushforward}. The degree-0 piece of this isomorphism induces
\[\iota':F_{k*}(T_0)=F_*(T)_0=\Int(F_*(T))_0\cong S_0.\]
By examining the diagram (\ref{construction-delta}) and the construction of the map $\delta$ in the exact sequence (\ref{L-ij-0-1}) in the cases $L=R/I$ and $L=\Int(F_*(R/I))$ respectively, the commutativity of the the middle square of (\ref{reduced-diagram}) can be reduced to the commutativity of the following diagram
\begin{equation}
\label{reduced-reduced-diagram}
\xymatrix{
F_{k*}T_0\ar[r]^{\mu_1}\ar[d]^{\iota'} & F_{k*}H^0(\mathbb{P}^n,\widetilde{T})\ar@{=}[r]\ar[d]^{\iota^{\prime\prime}}& F_{k*}H^0(\mathbb{P}^n,\underline{\Ext}^{1-j}_{\mathbb{P}^n}(\eta_*\cO_X,\dualp))\\
S_0\ar[r]^{\mu_2} & H^0(\mathbb{P}^n,\widetilde{S})\ar@{=}[r] & H^0(\mathbb{P}^n,\underline{\Ext}^{1-j}_{\mathbb{P}^n}(F_{\mathbb{P}^n*}\eta_*\cO_X,\dualp))
}
\end{equation}
where $\mu_1$ and $\mu_2$ are induced by the natural map $M_0\to H^0(\mathbb{P}^n,\widetilde{M})$ for any graded $R$-module $M$ (see (LCSC) for details); the map $\iota^{\prime\prime}$ is 
\[F_{k*}H^0(\mathbb{P}^n,\widetilde{T})=H^0(\mathbb{P}^n,F_{\mathbb{P}^n*}\widetilde{T})=H^0(\mathbb{P}^n,\widetilde{\Int(F_*T)})\cong H^0(\mathbb{P}^n,\widetilde{S}).\]
But the commutativity of diagram (\ref{reduced-reduced-diagram}) is clear, since the map $M_0\to H^0(\mathbb{P}^n,\widetilde{M})$ is functorial (see (LCSC) for details) and hence the degree-preserving isomorphism $\Int(F_*(T))\cong S$ induces a commutative diagram which is exactly diagram (\ref{reduced-reduced-diagram}). This proves that the square in the middle is commutative.\par 
It remains to prove that commutativity of the bottom square. Let $M$ denote $\Ext^{n+1-j}_R(\Int(F_*(R/I)),\Omega)$ and $N$ denote $\Ext^{n+1-j}_R(R/I,\Omega)$. From the proof of Proposition \ref{M-ij}, especially diagram (\ref{construction-delta}), we see that the vertical map\newline $H^j_{\fm}(\Int(F_*(R/I)))_0\to \E^{0,-j}(\Int(F_*(R/I)))_0$ in the bottom square is given by
\begin{align}
H^j_{\fm}(\Int(F_*(R/I)))_0\cong \D(M_0)&\to \D(H^0_{\fm}(M)_0)\notag\\
&\cong\Ext^{n+1}_R(M,\Omega)_0= \E^{0,-j}(\Int(F_*(R/I)))_0\notag
\end{align}
and the other vertical map $F_{k*}H^j_{\fm}(R/I)_0\to F_{k*}\E^{0,-j}(R/I)_0$ in the bottom square is given by
\begin{align}
F_{k*}H^j_{\fm}(R/I)_0\cong F_{k*}\D(N_0)&\to F_{k*}\D(H^0_{\fm}(N)_0)\notag\\
&\cong F_{k*}\Ext^{n+1}_R(N,\Omega)_0= F_{k*}\E^{0,-j}(R/I)_0.\notag
\end{align}
Given these two vertical maps and the horizontal maps $\rho$ and $\zeta_3$, the commutativity of the bottom square is rather straightforward to check.\par
This finishes the proof of the proposition.
\end{proof}

Finally we are ready to prove Theorem \ref{indp-f-M}.

\begin{proof}[Proof of Theorem \ref{indp-f-M}]
By virtue of Proposition \ref{frobenius-M-i-j}, we will consider $\varphi'$, which is the composition of $\tau_1$, $\tau_2$ and $\tau_3$ as in Discussion \ref{graded-p-structures-on-cE}. It is equivalent to proving that the degree-0 piece of each of $\tau_1$, $\tau_2$ and $\tau_3$ is independent of the embedding.\par
Since $\tau_3:F_*(\cE^{i,j}(R/I))\to \cE^{i,j}(R/I)$ is the identity map on the underlying abelian groups and from Proposition \ref{M-ij} we know that $\cE^{i,j}(R/I)_0$ is independent of the embedding, the degree-0 piece of $\tau_3$ is independent of the embedding. It follows from Proposition \ref{independent-tau-2} that the degree-0 piece of $\tau_2$ is also independent of the embedding. It remains to show that the degree-0 piece of $\tau_1$ is independent of the embedding. $\tau_1$ is induced by $R/I\to F_*(R/I)$ whose image is contained in $\Int(F_*(R/I))$. Considering $R/I\to L=\Int(F_*(R/I))$, one sees immediately from Proposition \ref{M-ij} that the degree-0 piece of $\tau_1$ is independent of the embedding. This completes the proof. 
\end{proof}

\section{Proof of Theorem \ref{stable-M-i-j}}
In this section we prove Theorem \ref{stable-M-i-j} which says that
\[\lambda_{i,j}(A)=\dim_k((\cE^{i,j}(R/I)_0)_s).\] 
Since $(\cE^{i,j}(R/I)_0)_s$ is defined in terms of the $p$-linear structure $\varphi_0$ on $\cE^{i,j}(R/I)_0$ and we have already proven that both $\cE^{i,j}(R/I)_0$ and
$\varphi_0$ do not depend on the embeddings, this shows that
$\lambda_{i,j}(A)$ depends only on $X$ and completes the proof of
our main result.\par 

Since any field extension will not affect Lyubeznik numbers (i.e.,
if $I\subset R=k[x_0,\dots,x_n]$ is an ideal, $K$ is a field extension
of $k$ and $A_K=(R/IR)_{(x_0,\dots,x_n)}$, then $\lambda_{ij}(A)= \lambda_{ij}(A_K)$), we may and we do assume that
the underlying field $k$ is algebraically closed.\par

Let $S$ be a subset of a $R$-module $M$, throughout this section
$R\langle S\rangle$ will denote the $R$-submodule of $M$ generated
by all the elements of $S$.\par 

Let $M$ be a graded $R$-module. If $M$ admits a graded $p$-linear structure $f$, then it is clear that
\[M_s:= \bigcap^{\infty}_{i=1}f^i(M)=\bigcap^{\infty}_{i=1}f^i_0(M_0)=(M_0)_s.\]
If the degree-0 piece $M_0$ of $M$ is a finite $k$-space, so is $M_s$, and hence $N=R\langle M_s\rangle$ is finitely generated. \par
As in Discussion \ref{equivalence}, a graded $p$-linear structure $f$ on $M$ induces a degree-preserving $R$-module homomorphism
$$\beta:F^*M\to M$$ given by $\beta(r\otimes m)=rf(m)$. There
results an inverse system
$$M\xleftarrow{\beta} F^*M\xleftarrow{F^*\beta}F^{2*}(M)\leftarrow\cdots\leftarrow F^{i*}(M)\xleftarrow{F^{i*}\beta}F^{(i+1)*}(M)\leftarrow\cdots.$$

\begin{rem}
\label{compatibility}
{\rm When $M=\cE^{i,j}(R/I)$, the $p$-linear structure $\varphi$ on $\cE^{i,j}(R/I)$ as in Discussion \ref{graded-p-structures-on-cE} induces, as in previous paragraph, an inverse system $\{F^{e*}(\cE^{i,j}(R/I))\}_e$ in which all transition maps are degree-preserving. There is also a second inverse system $\{\cE^{i,j}(R/I^{[p^e]})\}_e$ where the degree-preserving transition maps are induced by the natural surjection $R/I^{[p^{e+1}]}\to R/I^{[p^e]}$. It is rather straightforward to check that the natural degree-preserving isomorphisms
\[F^{e*}(\cE^{i,j}(R/I))\xrightarrow{\sim}\cE^{i,j}(R/I^{[p^e]})\]
provide an isomorphism between the inverse systems 
\[ \{F^{e*}(\cE^{i,j}(R/I))\}_e\ {\rm and}\ \{\cE^{i,j}(R/I^{[p^e]})\}_e.\] 
In particular, we have }
\[\sideset{^*}{}\varprojlim_eF^{e*}(\cE^{i,j}(R/I))\cong \sideset{^*}{}\varprojlim_e\cE^{i,j}(R/I^{[p^e]}).\]
\end{rem}

\begin{lem}
\label{limit} 
Let $M$ be a finite graded $R$-module with a graded $p$-linear
structure given by a degree-preserving $R$-linear map $\beta:F^*M\to M$. Let
$N=R\langle M_s\rangle$. Then the graded $p$-linear structure on $M$ restricts to a graded $p$-linear structure on $N$, and the inclusion $N\subseteq M$ induces
$$\sideset{^*}{}\varprojlim_{i\to\infty} F^{i*}N=\sideset{^*}{}\varprojlim_{i\to\infty}F^{i*}M.$$
\end{lem}
\begin{proof}[Proof]
Define $f:M\to M$ by $f(m)=\beta(1\otimes m)$, then by Discussion \ref{equivalence} we know $f$ is a graded $p$-linear structure on $M$. If $\sum_ja_jn_j\in N$ where $a_j\in R$ and $n_j\in M_s$, then
\[f(\sum_ja_jn_j)=\sum_ja^p_jf(n_j)\in N\] since $n_j\in M_s$. Hence
$f|_N$ is indeed a graded $p$-linear structure on $N$. \par 
Since $F^*$ is exact, the exact sequence
$$0\to N\to M\to M/N\to 0$$
induces an exact sequence of inverse systems
{\tiny
$$\begin{CD}
@.   @AAA     @AAA    @AAA\\
 0   @>>>   (F^{i*}N)_n @>>>  (F^{i*}M)_n @>F^{i*}(\pi)_n>>  (F^{i*}(M/N))_n @>>> 0 \\
@.  @AA(F^{i*}\beta|_{F^*N})_nA   @AA(F^{i*}\beta)_nA   @AA(F^{i*}\bar{\beta})_nA\\
 0  @>>>   (F^{(i+1)*}N)_n  @>>>   (F^{(i+1)*}M)_n @>>F_R^{i+1}(\pi)_n>  (F^{(i+1)*}(M/N))_n @>>> 0\\
 @.  @AAA   @AAA   @AAA\\
\end{CD} $$}
Since $M$ is a finitely generated graded $R$-module, $M_0$ is a finite $k$-space, hence $M_s$ is also a finite $k$-space. Since $N=R\langle M_s\rangle$ is finitely generated, so are $F^{i*}N$, and hence each $(F^{i*}N)_n$ is a finite $k$-space, so the inverse system
$\{(F^{i*}N)_n\}$ satisfies the Mittag-Leffler condition (page 191 in \cite{ag}), we have (by Proposition 9.1 on page 192 in \cite{ag}) an exact sequence
$$0\to \varprojlim (F^{i*}N)_n\to \varprojlim (F^{i*}M)_n\to \varprojlim (F^{i*}(M/N))_n\to 0$$
for each $n$, which induces another exact sequence
\begin{equation}
\label{ex-seq}
0\to \sideset{^*}{}\varprojlim F^{i*}N\to \sideset{^*}{}\varprojlim F^{i*}M\to \sideset{^*}{}\varprojlim F^{i*}(M/N)\to 0.
\end{equation}

We claim that there exists an integer $l$ such that $F^{i*}(M/N)\to\cdots\to M/N$ (or briefly, $F^{i*}(M/N)\to M/N$) has 0 image in degrees $\leq 1$ for all $i\geq l$ and we reason as follows. The homomorphism $F^{i*}(M/N)\to M/N$ is given by $r\otimes \bar{m}\mapsto \overline{rf^i(m)}$ for $r\in R^{(i)}\cong R$ and $m\in M$, thus,
$$\im(F^{i*}(M/N)\to M/N)=R\langle f^i(M/N)\rangle.$$
Since $M/N$ is finitely generated, there exists an integer $l_0$ such that $f^i(\bar{m})=0$ if $\deg(\bar{m})<0$ and $i\geq l_0$. This tells us that $$R\langle f^i(M/N)\rangle=R\langle f^i((M/N)_{\geq 0})\rangle=R\langle f^i((M/N)_0)\rangle+R\langle f^i((M/N)_{>0})\rangle,\ {\rm for}\ i\geq l_0.$$
Since $\deg(f(m))=p\deg(m)$, we get that $R\langle f^i((M/N)_{>0})\rangle\subset (M/N)_{\geq p^i}$ which does not have any element sitting in $(M/N)_{\leq 1}$. Since $k$ is algebraically closed, we have a descending chain of $k$-subspaces
$$f(M_0)\supset f^2(M_0)\supset\cdots\supset f^i(M_0)\supset f^{i+1}(M_0)\supset \cdots$$
Since $M_0$ is a finite $k$-space, this descending chain stabilizes, i.e. there exists an integer $t$ such that
$$f^t(M_0)=f^{t+1}(M_0)=\cdots.$$
Hence $N=R\langle f^t(M_0)\rangle$. Thus, if $i\geq t$, then $f^i((M/N)_0)=0$. Therefore, we can pick $l=\max\{l_o,t\}$. This finishes the proof of our claim.\par
Followed directly from our claim, the map $F^{(i+j)*}(M/N)\to\cdots\to F^{j*}(M/N)$ has 0 image in degrees $\leq p^j$ for all $j$ and $i\geq l$. Since the maps in the inverse system are degree-preserving, $\varprojlim_i(F^{i*}(M/N))_n=0$ for all $n$, hence $\sideset{^*}{_i}\varprojlim(F^{i*}(M/N))=0$. By the short exact sequence (\ref{ex-seq}), one has
$$\sideset{^*}{_i}\varprojlim F^{i*}M=\sideset{^*}{_i}\varprojlim F^{i*}N.$$
 \end{proof}

We need the following result from \cite{lcd}.
\begin{lem}{(Lemma 1.14 in \cite{lcd})}
\label{basis}
Let $k$ be an algebraically closed field of characteristic $p>0$, and let $V$ be
a finite dimensional $k$-vector space with a bijective $p$-linear endomorphism
$f$. Then there is a basis $e_1,\cdots,e_r$ of $V$ with $f(e_i)=e_i$ for each $i$.
\end{lem}

\begin{pro}
\label{limit-stable}
Let $M,N$ be as in Lemma \ref{limit}. Let $M_s$ denote the stable part of $M$. Then
$$\sideset{^*}{}\varprojlim_iF^{i*}N=R^s,$$ where
$s=\dim_k(M_s)$.
\end{pro}
\begin{proof}[Proof]
Use induction on $s$. Notice that $N$ is generated by $s$ element and clearly $f$ acts bijectively on $M_s$.\par
When $s=1$, $N$ is generated by a single element $a\in M_0$ with $f(a)=a$ by Lemma \ref{basis}. Then $N\cong R/I$,
where $I=\Ann_R(a)$. Hence we have
$$\sideset{^*}{}\varprojlim_{i\to\infty}F^{i*}N\cong\sideset{^*}{}\varprojlim_{i\to\infty}R/I^{[p^i]}= R,$$
where the last equality follows from Example \ref{example}.\par
Assume that $s\geq 2$ and the proposition is true for $s-1$. Let
$a_1,\cdots,a_s$ denote the generators of $N$ with $f(a_i)=a_i$. They exist by Lemma \ref{basis}. Let $N'$ be
the submodule of $N$ generated by $a_s$. Then, the short exact sequence
\[0\to N'\to N\to N/N'\to 0\]
induces, as in the proof of Lemma \ref{limit}, an exact sequence
\[0\to \sideset{^*}{}\varprojlim F^{i*}N'\to \sideset{^*}{}\varprojlim
F^{i*}N\to \sideset{^*}{}\varprojlim F^{i*}(N/N')\to 0.\]
is an exact sequence. By the inductive hypothesis, $\sideset{^*}{}\varprojlim
F^{i*}(N/N')\cong R^{s-1}$ which is projective over
$R$, and hence the exact sequence splits,
$$\sideset{^*}{}\varprojlim F^{i*}N\cong \sideset{^*}{}\varprojlim
F^{i*}N'\oplus R^{s-1}\cong R\oplus R^{s-1}\cong
R^s.$$
\end{proof}

\begin{pro}
\label{limit-lambda}
$\sideset{^*}{}\varprojlim_eF^{e*}(\cE^{i,j}(R/I))\cong R^{\lambda_{i,j}(A)}$.
\end{pro}
\begin{proof}[Proof]
Let $\E$ denote the injective hull of $R/\mathfrak{m}$ and $\sideset{^*}{}\E$ denote the $^*$injective hull of $R/\mathfrak{m}$ (cf. 13.2.1 in \cite{bs}). First we notice that
\begin{align}
 &\sideset{^*}{_R}\Hom(H^i_{\mathfrak{m}}(H^{n+1-j}_I(R))(-n-1),\sideset{^*}{}\E)\notag\\
                           &\cong \sideset{^*}{_R}\Hom(H^i_{\mathfrak{m}}(\varinjlim_e\Ext^{n+1-j}_R(R/I^{[p^e]},R))(-n-1),\sideset{^*}{}\E)\tag{i}\\
                           &\cong \sideset{^*}{_R}\Hom(\varinjlim_eH^i_{\mathfrak{m}}(\Ext^{n+1-j}_R(R/I^{[p^e]},\Omega)),\sideset{^*}{}\E)\tag{ii}\\
                           &\cong \sideset{^*}{_e}\varprojlim\sideset{^*}{_R}\Hom(H^i_{\mathfrak{m}}(\Ext^{n+1-j}_R(R/I^{[p^e]},\Omega))),\sideset{^*}{}\E)\tag{iii}\\
                          &\cong \sideset{^*}{_e}\varprojlim\cE^{i,j}(R/I^{[p^e]})\tag{iv}\\
                          &\cong \sideset{^*}{_e}\varprojlim(F^{e*}(\cE^{i,j}(R/I)))\tag{v}
\end{align}
(i) follows from the definition of local cohomology.\\
(ii) holds since local cohomology commutes with direct limits and degree shifting.\\
(iii) follows from Proposition \ref{graded-inverse}.\\
(iv) is a consequence of (GLD).\\
(v) is nothing but Remark \ref{compatibility}.\\
By Proposition \ref{limit-stable}, we see that $\sideset{^*}{_e}\varprojlim(F^{e*}(\cE^{i,j}(R/I)))=R^s$, as graded $R$-modules, where $s$ is the dimension of the stable part of $\cE^{i,j}(R/I)$. On the other hand, according to Lemma 2.2 in \cite{l3}, $H^i_{\mathfrak{m}}(H^{n+1-j}_I(R))\cong \E^{\lambda_{i,j}(A)}$ (regardless of the grading). It follows immediately that 
\[H^i_{\mathfrak{m}}(H^{n+1-j}_I(R))(-n-1)\cong \sideset{^*}{^{\lambda_{i,j}(A)}}\E.\]
Therefore 
\begin{align}
R^{\lambda_{i,j}(A)}\cong  \sideset{^*}{_R}\Hom(\sideset{^*}{^{\lambda_{i,j}(A)}}\E,\sideset{^*}{}\E)&\cong \sideset{^*}{_R}\Hom(H^i_{\mathfrak{m}}(H^{n+1-j}_I(R))(-n-1),\sideset{^*}{}\E)\notag\\
&\cong \sideset{^*}{_e}\varprojlim(F^{e*}(\cE^{i,j}(R/I))).\notag
\end{align}
\end{proof}

\begin{proof}[Proof of Theorem \ref{stable-M-i-j}]
By Proposition \ref{limit-stable} and Proposition \ref{limit-lambda}, we have
$$\lambda_{i,j}(A)=\dim_k(\cE^{i,j}(R/I)_s)=\dim_k(\cE^{i,j}(R/I)_0)_s).$$
\end{proof}

\section*{Acknowledgments}
This paper develops and completes the results in my old preprint \cite{zh1} which was a part of my thesis. I would like to thank my advisor, Professor Gennady Lyubeznik, for his support, guidance and numerous discussions of the material in this paper. I am grateful to Professor Joseph Lipman for answering my questions on dualizing complexes and duality, for carefully reading an earlier draft of this paper, and for pointing out a generalization of Lemma \ref{dualizing-structure}. Special thanks are due to Karl Schwede for carefully reading an earlier draft of this paper and for his helpful suggestions. I also want to thank the referee for his/her comments that improve the exposition of this paper considerably.


\end{document}